\documentclass{amsart}

\usepackage{amsthm}
\usepackage[leqno]{amsmath}
\usepackage{latexsym,amsfonts,amssymb}
\usepackage[all]{xy} \SelectTips{eu}{} \SilentMatrices
\usepackage{hyperref}

\newcommand{\numberseries}{\mdseries}   

\newlength{\thmtopspace}                
\newlength{\thmbotspace}                
\newlength{\thmheadspace}               
\newlength{\thmindent}                  

\renewcommand{\subparagraph}{\vspace{\thmbotspace}}

\newtheoremstyle{bfupright head,slanted body}
                {\thmtopspace}{\thmbotspace}
                {\slshape}{\thmindent}{\bfseries}{.}{\thmheadspace}
                {{\numberseries \thmnumber{(#2) }}\thmnote{#3}}

\newtheoremstyle{bfupright head,upright body}
                {\thmtopspace}{\thmbotspace}
                {\upshape}{\thmindent}{\bfseries}{.}{\thmheadspace}
                {{\numberseries \thmnumber{(#2) }}\thmnote{#3}}

\newtheoremstyle{fixed bf head,slanted body}
                {\thmtopspace}{\thmbotspace}{\slshape}
                {\thmindent}{\bfseries}{.}{\thmheadspace}
                {{\numberseries \thmnumber{(#2) }}\thmname{#1}\thmnote{ (#3)}}

\newtheoremstyle{fixed bf head,upright body}
                {\thmtopspace}{\thmbotspace}{\upshape}
                {\thmindent}{\bfseries}{.}{\thmheadspace}
                {{\numberseries \thmnumber{(#2) }}\thmname{#1}\thmnote{ (#3)}}

\newtheoremstyle{numbered paragraph}
                {\thmtopspace}{\thmbotspace}{\upshape}
                {\thmindent}{\upshape}{}{0pt}
                {{\numberseries \thmnumber{(#2) }}}

\newtheoremstyle{unnumbered paragraph}
                {\thmtopspace}{\thmbotspace}{\upshape}
                {\parindent}{\upshape}{}{0pt}

\theoremstyle{bfupright head,slanted body}
\newtheorem{res}{}[section]             \newtheorem*{res*}{}
\theoremstyle{bfupright head,upright body}
\newtheorem{bfhpg}[res]{}               \newtheorem*{bfhpg*}{}
\theoremstyle{fixed bf head,slanted body}
\newtheorem{thm}[res]{Theorem}          \newtheorem*{thm*}{Theorem}
\newtheorem{prp}[res]{Proposition}      \newtheorem*{prp*}{Proposition}
\newtheorem{cor}[res]{Corollary}        \newtheorem*{cor*}{Corollary}
\newtheorem{lem}[res]{Lemma}            \newtheorem*{lem*}{Lemma}

\theoremstyle{fixed bf head,upright body}
\newtheorem{dfn}[res]{Definition}       \newtheorem*{dfn*}{Definition}
\newtheorem{con}[res]{Construction}     \newtheorem*{con*}{Construction}
      \newtheorem*{obs*}{Observation}
\newtheorem{rmk}[res]{Remark}           \newtheorem*{rmk*}{Remark}
\newtheorem{exa}[res]{Example}          \newtheorem*{exa*}{Example}
\theoremstyle{numbered paragraph}
\newtheorem{ipg}[res]{}

\newlength{\thmlistleft}        
\newlength{\thmlistright}       
\newlength{\thmlistpartopsep}   
\newlength{\thmlisttopsep}      
\newlength{\thmlistparsep}      
\newlength{\thmlistitemsep}     

\newcounter{eqc} 
\newenvironment{eqc}{\begin{list}{\upshape (\textit{\roman{eqc}})}%
    {\usecounter{eqc}%
      \setlength{\leftmargin}{\thmlistleft}%
      \setlength{\labelwidth}{\thmlistleft}%
      \setlength{\rightmargin}{\thmlistright}%
      \setlength{\partopsep}{\thmlistpartopsep}%
      \setlength{\topsep}{\thmlisttopsep}%
      \setlength{\parsep}{\thmlistparsep}%
      \setlength{\itemsep}{\thmlistitemsep}}}%
  {\end{list}}%

\newcommand{\eqclbl}[1]{{\upshape(\textit{#1})}}

\newcounter{prt}
  {\end{list}}%

\newenvironment{prf*}[1][Proof]{%
  \begin{proof}[\bf #1]
    \setcounter{equation}{0}
    }
  {\end{proof}
}
  \newcommand{\proofofimp}[3][:]{\mbox{\eqclbl{#2}$\!\implies\!$\eqclbl{#3}#1}}

\newcommand{\pgref}[1]{(\ref{#1})}
\newcommand{\thmref}[2][Theorem~]{#1\pgref{thm:#2}}
\newcommand{\corref}[2][Corollary~]{#1\pgref{cor:#2}}
\newcommand{\prpref}[2][Proposition~]{#1\pgref{prp:#2}}
\newcommand{\lemref}[2][Lemma~]{#1\pgref{lem:#2}}

\newcommand{\conref}[2][Construction~]{#1\pgref{con:#2}}
\newcommand{\dfnref}[2][Definition~]{#1\pgref{dfn:#2}}

\newcommand{\rmkref}[2][Remark~]{#1\pgref{rmk:#2}}
\renewcommand{\eqref}[1]{\pgref{eq:#1}}
\newcommand{\thmcite}[2][?]{\cite[thm.~#1]{#2}}
\newcommand{\prpcite}[2][?]{\cite[prop.~#1]{#2}}
\newcommand{\lemcite}[2][?]{\cite[lem.~#1]{#2}}
\newcommand{\seccite}[2][?]{\cite[sec.~#1]{#2}}
\newcommand{\set}[2][\mspace{1mu}]{\{#1 #2 #1\}}
\newcommand{\setof}[3][\mspace{1mu}]{\{#1#2 \mid #3#1\}}
\newcommand{\ZZ}{\mathbb{Z}}
\newcommand{\qtext}[1]{\quad\text{#1}\quad}
\newcommand{\qqtext}[1]{\qquad\text{#1}\qquad}
\newcommand{\qand}{\qtext{and}}
\newcommand{\qqand}{\qqtext{and}}
\newcommand{\deq}{\:=\:}
\DeclareMathOperator*{\dfinsum}{\textstyle\bigoplus}
\renewcommand{\a}{\alpha}
\renewcommand{\b}{\beta}
\newcommand{\f}{\varphi}
\newcommand{\m}{\mathfrak{m}}
\newcommand{\is}{\cong}
\newcommand{\qis}{\simeq}
\renewcommand{\le}{\leqslant}
\renewcommand{\ge}{\geqslant}

\newcommand{\lra}{\longrightarrow}
\newcommand{\xra}[2][]{\xrightarrow[#1]{\;#2\;}}
\newcommand{\qra}{\xra{\;\qis\;}}
\newcommand{\poly}[2][k]{#1[#2]}
\newcommand{\mapdef}[4][\rightarrow]{\nobreak{#2\colon #3 #1 #4}}
\newcommand{\dmapdef}[4][\lra]{\nobreak{#2\colon #3\:#1\:#4}}
\newcommand{\Coker}[1]{\nobreak{\operatorname{Coker}#1}}
\newcommand{\Cone}[1]{\nobreak{\operatorname{Cone}#1}}
\newcommand{\dgr}[1]{|#1|}
\newcommand{\dif}[2][]{{\partial}^{#2}_{#1}}
\renewcommand{\Im}[1]{\operatorname{Im}{#1}}
\newcommand{\Bo}[2][]{\operatorname{B}_{#1}(#2)}
\newcommand{\Cy}[2][]{\operatorname{Z}_{#1}(#2)}
\newcommand{\Co}[2][]{\operatorname{C}_{#1}(#2)}
\renewcommand{\H}[2][]{\operatorname{H}_{#1}(#2)}
\newcommand{\Shift}[2][]{\mathsf{\Sigma}^{#1}{#2}}
\newcommand{\Shiftp}[2][]{(\Shift[#1]{#2})}
\newcommand{\Tha}[2]{#2_{{\scriptscriptstyle\le}#1}}
\newcommand{\Thb}[2]{#2_{{\scriptscriptstyle\ge}#1}}
\newcommand{\Tsa}[2]{#2_{{\scriptscriptstyle\subset}#1}}
\newcommand{\Tsb}[2]{#2_{{\scriptscriptstyle\supset}#1}}
\newcommand{\Thap}[2]{(\Tha{#1}{#2})}
\newcommand{\Thbp}[2]{(\Thb{#1}{#2})}

\newcommand{\rnk}[2][k]{\operatorname{rank}_{#1}#2}
\newcommand{\fd}[2][R]{\operatorname{fd}_{#1}#2}
\newcommand{\id}[2][R]{\operatorname{id}_{#1}#2}
\newcommand{\pd}[2][R]{\operatorname{pd}_{#1}#2}
\newcommand{\Gid}[2][R]{\operatorname{Gid}_{#1}#2}
\newcommand{\Gpd}[2][R]{\operatorname{Gpd}_{#1}#2}

\newcommand{\Hom}[3][R]{\operatorname{Hom}_{#1}(#2,#3)}
\newcommand{\RHom}[3][R]{\operatorname{\mathbf{R}Hom}_{#1}(#2,#3)}
\newcommand{\Ext}[4][R]{\operatorname{Ext}_{#1}^{#2}(#3,#4)}
\newcommand{\tp}[3][R]{\nobreak{#2\otimes_{#1}#3}}

\newcommand{\tpp}[3][R]{(\tp[#1]{#2}{#3})}
\newcommand{\Ltp}[3][R]{\nobreak{#2\otimes_{#1}^{\mathbf{L}}#3}}
\newcommand{\Tor}[4][R]{\operatorname{Tor}^{#1}_{#2}(#3,#4)}
\newcommand{\Cat}[2]{{\mathsf{#2}}(#1)}
\newcommand{\C}[1][R]{\Cat{#1}{C}}
\newcommand{\D}[1][R]{\Cat{#1}{D}}

\numberwithin{equation}{res}

\newcommand{\bicat}[2]{#1\text{--}#2^{\circ}}
\newcommand{\bi}[3][bi]{\mbox{$#2$--$#3^{\circ}$}\nobreakdash-#1}
\def\urltilda{\kern -.15em\lower .7ex\hbox{\~{}}\kern .04em}
\newcommand{\btp}[3][R]{\nobreak{#2\otimes^\Join_{#1}#3}}
\newcommand{\pHom}[3][R]{\operatorname{Hom}^\Join_{#1}(#2,#3)}

\newcommand{\btpp}[3][R]{\nobreak{(#2\otimes^\Join_{#1}#3)}}
\newcommand{\Ttor}[4][R]{\smash{\operatorname{\widehat{Tor}}}_{#2}^{#1{^{\phantom{|\mspace{-6mu}}}
    }}(#3,#4)}
\newcommand{\Text}[4][R]{\smash{\operatorname{\widehat{Ext}}}_{#1}^{#2{^{\phantom{|}\mspace{-6mu}}}}(#3,#4)}

\newcommand{\Rop}{R^{\circ}}
\newcommand{\Sop}{S^{\circ}}

\def\widebardisplay#1{%
  \setbox0=\hbox{$\displaystyle #1$}
  \dimen0=\wd0%
  \advance\dimen0 by -2pt
  \vbox{%
    \nointerlineskip%
    \moveright 2pt 
    \vbox{\hrule width \dimen0}%
    \nointerlineskip%
    \kern 1pt
    \box0%
    }%
  }

\def\widebartext#1{%
  \setbox0=\hbox{$#1$}
  \dimen0=\wd0%
  \advance\dimen0 by -2pt
  \vbox{%
    \nointerlineskip%
    \moveright 1pt 
    \vbox{\hrule width \dimen0}%
    \nointerlineskip%
    \kern 1.2pt
    \box0%
    }%
  }

\def\widebarscript#1{%
  \setbox0=\hbox{$\scriptstyle #1$}
  \dimen0=\wd0%
  \advance\dimen0 by -1.5pt
  \vbox{%
    \nointerlineskip%
    \moveright .75pt 
    \vbox{\hrule width \dimen0}%
    \nointerlineskip%
    \kern .9pt
    \box0%
    }%
  }

\def\widebarscriptscript#1{%
  \setbox0=\hbox{$\scriptscriptstyle #1$}
  \dimen0=\wd0%
  \advance\dimen0 by -2pt
  \vbox{%
    \nointerlineskip%
    \moveright 1pt 
    \vbox{\hrule width \dimen0}%
    \nointerlineskip%
    \kern .6pt
    \box0%
    }%
  }

\def\widebar#1{\mathchoice%
  {\widebardisplay{#1}}%
  {\widebartext{#1}}%
  {\widebarscript{#1}}%
  {\widebarscriptscript{#1}}%
  }

\makeatletter
\renewcommand{\intertext}{\@amsmath@err{\Invalid@@\intertext}\@eha}
\def\intertext@{%
  \def\intertext##1{%
    \ifvmode\else\\\@empty\fi
    \noalign{%
      \penalty\postdisplaypenalty\vskip0.3\belowdisplayskip
      \vbox{\normalbaselines
        \ifdim\linewidth=\columnwidth
        \else \parshape\@ne \@totalleftmargin \linewidth
        \fi
        \noindent##1\par}%
      \penalty\predisplaypenalty\vskip0.3\belowdisplayskip%
    }%
}}
\makeatother

\begin{document}

\title{Tate (co)homology via pinched complexes}

\author[Lars Winther Christensen]{Lars Winther Christensen}

\address{Department of Math.\ and Stat., Texas Tech University,
  Lubbock, TX 79409, U.S.A.}

\email{lars.w.christensen@ttu.edu}

\urladdr{http://www.math.ttu.edu/\urltilda lchriste}

\author[David A.~Jorgensen]{David A. Jorgensen}

\address{Department of Mathematics, University of Texas, Arlington,
  TX~76019, U.S.A.}

\email{djorgens@uta.edu}

\urladdr{http://dreadnought.uta.edu/\urltilda dave}

\thanks{Work partly supported by NSA grants H98230-11-0214 (L.W.C.)
  and H98230-10-0197~(D.A.J.).}

\date{13 November 2011}





\begin{abstract}
  For complexes of modules we study two new constructions, which we
  call the pinched tensor product and the pinched Hom. They provide
  new methods for computing Tate homology
  $\smash{\operatorname{\widehat{Tor}}}$ and Tate cohomology
  $\smash{\operatorname{\widehat{Ext}}}$, which lead to conceptual
  proofs of balancedness of Tate (co)homology for modules over
  associative~rings.

  Another application we consider is in local algebra. Under
  conditions of vanishing of Tate (co)homology, the pinched tensor
  product of two minimal complete resolutions yields a minimal
  complete resolution.
\end{abstract}

\maketitle


\thispagestyle{empty}

\section*{Introduction}
\label{sec:intro}

\noindent
Tate cohomology originated in the study of representations of finite
groups. It has been generalized---through works of, in chronological
order, Buchweitz~\cite{ROB86}, Avramov and Martsinkovsky
\cite{LLAAMr02}, and Veliche \cite{OVl06}---into a cohomology theory
for modules with complete resolutions. The parallel theory of Tate
homology has been treated in the same generality by Iacob
\cite{AIc07}.

While these theories function for modules over any associative ring,
the central question of balancedness has yet to receive a cogent
treatment. The extant literature only solves the problem for modules
over special commutative rings. The issue is that if $M$ and $N$ are
modules with appropriate complete resolutions, then there are
potentially two ways of defining Tate cohomology $\Text[]{*}{M}{N}$;
do they yield the same theory? For Tate homology $\Ttor[]{*}{M}{N}$
one encounters a similar situation, and one goal of this paper is to
resolve these balancedness problems.

Proving balancedness of absolute (co)homology, Ext and Tor, boils down
to showing that, say, $\Tor[]{*}{M}{N}$ can be computed from a complex
constructed from resolutions of both variables $M$ and $N$; namely the
tensor product of their projective resolutions. Our approach is
similar, but for Tate (co)homology the standard tensor product and Hom
complexes fail to do the job, so we introduce two new
constructions. We call them the pinched tensor product and the pinched
Hom. They resemble the usual tensor product and Hom of complexes, but
they are smaller in a sense that is discussed below. The central
technical results are \thmref[Theorems~]{btp-Ttor} and
\thmref[]{pHom-Text1}, which establish that Tate (co)homology can be
computed from pinched complexes. The balancedness problems are
resolved in \thmref[Theorems~]{Ttor-bal} and \thmref[]{Text-bal}.

As part of our analysis of the pinched complexes, we establish
``pinched versions'' of standard isomorphisms for complexes, such as
Hom-tensor adjunction. They allow us to give
criteria---\corref[Corollaries~]{pp=p} and \corref[]{pi=i}---in terms
of vanishing of Tate (co)homology, for when a pinched Hom complex
$\pHom[]{T}{U}$ or a pinched tensor product $\btp[]{T}{U}$ of complete
resolutions is a complete resolution.

This is of particular interest in local algebra since, if one starts
with unbounded complexes of finitely generated modules, then the
pinched Hom and the pinched tensor product are also complexes of
finitely generated modules.  \thmref{TR} gives a criterion, in terms
of vanishing of Tate (co)homology, for a tensor product of minimal
complete resolutions to be a minimal complete resolution.

\section{Standard constructions with complexes}
\label{sec:thac}

\noindent
In this paper $R$, $R'$, $S$, and $S'$ are associative unital rings;
they are assumed to be algebras over a common commutative unital ring
$\Bbbk$.  The default $\Bbbk$ is the ring $\ZZ$ of integers, but in
concrete settings other choices may be useful. For example, if the
rings are algebras over a field $k$, then $\Bbbk = k$ is a natural
choice. If $R$ is commutative, and $R'$, $S$, and $S'$ are
$R$-algebras, then $\Bbbk = R$ is a candidate.

Modules are assumed to be unitary, and the default action of the ring
is on the left. Right modules over $R$ are hence treated as (left)
modules over the opposite ring $\Rop$. By an \bi{R}{S}module we mean a
module over the $\Bbbk$-algebra $\tp[\Bbbk]{R}{\Sop}$. Note that every
$R$-module has a natural \bi{R}{\Bbbk}module structure; in particular
they are symmetric \bi{\Bbbk}{\Bbbk}modules.  Modules over a
commutative ring $R$ are tacitly assumed to be symmetric
\bi{R}{R}modules.

\begin{bfhpg*}[Complexes]
  An $R$-complex is a (homologically) graded $R$-module $M$ endowed
  with a square-zero endomorphism $\dif{M}$ of degree $-1$, which is
  called the differential. Here is a visualization,
  \begin{equation*}
    \cdots \lra M_{i+1} \xra{\dif[i+1]{M}} M_i \xra{\dif[i]{M}}
    M_{i-1} \lra \cdots
  \end{equation*}
  A morphism of complexes $M \to N$ is a degree $0$ graded
  homomorphism $\a = (\a_i)_{i\in\ZZ}$ of the underlying graded
  modules that commutes with the differentials on $M$ and $N$; i.e.\
  one has $\dif{N}\a = \a\dif{M}$.  The category of $R$-complexes is
  denoted $\C$. If the underlying graded module is an \bi{R}{S}module,
  and the differential is a bimodule endomorphism, then the complex is
  called a complex of \bi{R}{S}modules;\footnote{ The term `bicomplex'
    is too close to `double complex'.} the category of such complexes
  is denoted $\C[\bicat{R}{S}]$.

  The kernel $\Cy{M}$ and the image $\Bo{M}$ of $\dif{M}$ are graded
  submodules of $M$ and, in fact, subcomplexes, as the induced
  differentials are trivial. A complex $M$ is called \emph{acyclic} if
  the homology complex $\H{M} = \Cy{M}/\Bo{M}$ is the zero-complex. We
  use the notation $\Co{M}$ for the cokernel of the differential,
  i.e.\ $\Co[i]{M}= \Coker{\dif[i+1]{M}}$.

  The notation $\sup{M}$ and $\inf{M}$ is used for the supremum and
  infimum of the set $\setof{i \in \ZZ}{M_i \ne 0}$, with the
  conventions $\sup{\emptyset} = -\infty$ and $\inf{\emptyset} =
  \infty$. A complex $M$ is \emph{bounded above} if $\sup{M}$ is
  finite, and it is \emph{bounded below} if $\inf{M}$ is finite.

  For $n\in\ZZ$ the \emph{n-fold shift} of $M$ is the complex
  $\Shift[n]{M}$ with $(\Shift[n]{M})_i = M_{i-n}$ and
  $\dif[i]{\Shift[n]{M}} = (-1)^n\dif[i-n]{M}$. One has
  $\sup{\Shiftp[n]{M}} = \sup{M}+n$ and $\inf{\Shiftp[n]{M}} =
  \inf{M}+n$.

  Let $n$ be an integer. The \emph{hard truncation above} of $M$ at
  $n$ is the complex $\Tha{n}{M}$ with $\Thap{n}{M}_i=0$ for $i >n$
  and $\dif[i]{\Tha{n}{M}} = \dif[i]{M}$ for $i \le n$. It looks like
  this:
  \begin{equation*}
    \Tha{n}{M} \deq 0 \lra M_n \xra{\dif[n]{M}} M_{n-1} \xra{\dif[n-1]{M}}
    M_{n-2}\lra \cdots
  \end{equation*}
  Similarly, $\Thb{n}{M}$ is the complex with $\Thbp{n}{M}_i=0$ for
  $i<n$ and $\dif[i]{\Thb{n}{M}} = \dif[i]{M}$ for~$i>n$. Note that
  $\Tha{n}{M}$ is a subcomplex of $M$, and $\Thb{n}{M}$ is the
  quotient complex $M/\Tha{n-1}{M}$. The \emph{soft truncations} of
  $M$ at $n$ are the complexes
  \begin{align*}
    \Tsa{n}{M} & \deq 0 \lra \Co[n]{M} \xra{\widebar{\dif[n]{M}}}
    M_{n-1} \xra{\dif[n-1]{M}} M_{n-2}\lra \cdots \intertext{and}
    \Tsb{n}{M} &\deq \cdots \lra M_{n+2} \xra{\dif[n+2]{M}} M_{n+1}
    \xra{\dif[n+1]{M}} \Cy[n]{M} \lra 0.
  \end{align*}

  A morphism of complexes that induces an isomorphism in homology is
  called a \emph{quasi-isomorphism} and indicated by the symbol
  `$\qis$'. A morphism $\a$ is a quasi-isomorphism if and only if its
  mapping cone, the complex $\Cone{\a}$, is acyclic.
\end{bfhpg*}

The central constructions in this paper, \conref[]{btp} and
\conref[]{pHom}, start from the standard constructions of tensor
product and Hom complexes, hence we review them in~detail.

\begin{bfhpg*}[Tensor~product~and~Hom]
  Let $M$ be an $\Rop$-complex and $N$ be an $R$-complex. The tensor
  product $\tp{M}{N}$ is the $\Bbbk$-complex whose underlying graded
  module is given by
  \begin{equation*}
    \tpp{M}{N}_n = \coprod_{i\in\ZZ} \tp{M_i}{N_{n-i}},
  \end{equation*}
  and whose differential is defined by specifying its action on an
  elementary tensor of homogeneous elements as follows,
  \begin{equation*}
    \dif{\tp{M}{N}}(x\otimes y)= \dif{M}(x)\otimes y + (-1)^{\dgr{x}} x\otimes
    \dif{N}(y);
  \end{equation*}
  here $\dgr{x}$ is the degree of $x$ in $M$. For a morphism of
  $\Rop$-complexes $\mapdef{\a}{M}{M'}$ and a morphism of
  $R$-complexes $\mapdef{\b}{N}{N'}$, the map
  $\mapdef{\tp{\a}{\b}}{\tp{M}{N}}{\tp{M'}{N'}}$, defined~by
  \begin{equation*}
    (\tp{\a}{\b})(x\otimes y) = \a(x)\otimes\b(y),
  \end{equation*}
  is a morphism of $\Bbbk$-complexes. The tensor product yields a
  functor
  \begin{equation*}
    \dmapdef{\tp{-}{-}}{\C[\Rop] \times
      \C}{\C[\Bbbk]},
  \end{equation*}
  which is $\Bbbk$-bilinear and right exact in each variable.  In case
  $M$ is a complex of \bi{R'}{R}modules and $N$ is a complex of
  \bi{R}{S}modules, then the tensor product $\tp{M}{N}$ is a complex
  of \bi{R'}{S}modules. The tensor product yields a functor
  $\C[\bicat{R'}{R}] \times \C[\bicat{R}{S}] \to \C[\bicat{R'}{S}]$.

  \subparagraph For $R$-complexes $M$ and $N$, the $\Bbbk$-complex
  $\Hom{M}{N}$ is given by
  \begin{align*}
    \Hom{M}{N}_n &= \prod_{i\in\ZZ}\Hom{M_i}{N_{i+n}} \intertext{and}
    \dif{\Hom{M}{N}}(\f) &= \dif{N}\f -(-1)^{\dgr{\f}}\f\dif{M},
  \end{align*}
  for a homogeneous $\f$ in $\Hom{M}{N}$. For morphisms of
  $R$-complexes $\mapdef{\a}{M}{M'}$ and $\mapdef{\b}{N}{N'}$, a
  morphism $\mapdef{\Hom{\a}{\b}}{\Hom{M'}{N}}{\Hom{M}{N'}}$ of
  \linebreak $\Bbbk$\nobreakdash-complexes is defined~by
  \begin{equation*}
    \Hom{\a}{\b}(\f)=\b\f\a.
  \end{equation*}
  With these definitions, Hom yields a functor,
  \begin{equation*}
    \dmapdef{\Hom{-}{-}}{{\C}^\mathrm{op} \times
      \C}{\C[\Bbbk]},
  \end{equation*}
  where the superscript `op' signifies the opposite category; it is
  $\Bbbk$-bilinear and left exact in each variable. In case $M$ is a
  complex of \bi{R}{R'}modules and $N$ is a complex of
  \bi{R}{S}modules, the complex $\Hom{M}{N}$ is one of
  \bi{R'}{S}modules; Hom yields a functor
  $\C[\bicat{R}{R'}]^\mathrm{op} \times \C[\bicat{R}{S}] \to
  \C[\bicat{R'}{S}]$.
\end{bfhpg*}

\begin{bfhpg*}[Resolutions]
  \label{res}
  An $R$-complex $P$ is called \emph{semi-projective} if each module
  $P_i$ is projective, and the functor $\Hom{P}{-}$ preserves
  quasi-isomorphisms (equivalently, it preserves acyclicity). A
  bounded below complex of projective $R$-modules is
  semi-projective. Similarly, an $R$-complex $I$ is called
  \emph{semi-injective} if each module $I_i$ is injective, and the
  functor $\Hom{-}{I}$ preserves quasi-isomorphisms (equivalently, it
  preserves acyclicity). A bounded above complex of injective
  $R$-modules is semi-injective. Every $R$-complex $M$ has a
  semi-projective resolution and a semi-injective resolution; that is,
  there are quasi-isomorphisms $\mapdef{\pi}{P}{M}$ and
  $\mapdef{\iota}{M}{I}$, where $P$ is semi-projective and $I$ is
  semi-injective; see \cite{LLAHBF91}\footnote{ In this paper the
    authors use 'DG-' in place of 'semi-'.}. An $R$-complex $F$ is
  called \emph{semi-flat} if each module $F_i$ is flat, and the
  functor $\tp{-}{F}$ preserves quasi-isomorphisms (equivalently, it
  preserves acyclicity). Every semi-projective complex is semi-flat.

  For an $R$-module $M$, a projective (injective) resolution in the
  classic sense is a semi-projective (-injective) resolution. Thus,
  the following definitions of homological dimensions of an
  $R$-complex extend the classic notions for modules.
  \begin{align*}
    \pd{M} &= \inf\setof{\sup{P}}{P \qra M \text{ is a semi-projective
        resolution}},\\
    \id{M} &= \inf\setof{-\inf{I}}{M \qra I \text{ is a semi-injective
        resolution}},\qand\\
    \fd{M} &= \inf\left\{\:n\ge \sup{\H{M}}\:\left|\:
        \begin{gathered}
          \text{ $P \qra M$ is a semi-projective} \\
          \text{ resolution and $\Co[n]{P}$ is flat }
        \end{gathered}
      \right.  \right\}.
  \end{align*}

  The derived tensor product $\Ltp{-}{-}$ and the derived Hom functor
  $\RHom{-}{-}$ for complexes are computed by way of the resolutions
  described above. Extending the usual definitions of Tor and Ext for
  modules, set
  \begin{equation*}
    \Tor{i}{M}{N} = \H[i]{\Ltp{M}{N}} \qand \Ext{i}{M}{N} =
    \H[-i]{\RHom{M}{N}}
  \end{equation*}
  for complexes $M$ and $N$ and $i\in\ZZ$.
\end{bfhpg*}

\section{Complete resolutions and Tate homology}
\label{sec:tatehomology}

\noindent
In this section we recall some definitions and facts from works of
Iacob~\cite{AIc07} and Veliche~\cite{OVl06}, and we establish some
auxiliary results for later use.

\begin{bfhpg}[Complete projective resolutions]
  An acyclic complex $T$ of projective $R$\nobreakdash-modules is
  called \emph{totally acyclic,} if the complex $\Hom{T}{Q}$ is
  acyclic for every projective $R$-module $Q$.

  A \emph{complete projective resolution} of an $R$-complex $M$ is a
  diagram
  \begin{equation}
    \label{eq:cpltres}
    T \xra{\tau} P \xra{\pi} M,
  \end{equation}
  where $\pi$ is a semi-projective resolution, $T$ is a totally
  acyclic complex of projective $R$-modules, and $\tau_i$ is an
  isomorphism for $i \gg 0$.
\end{bfhpg}

See \cite{OVl06} for a proof of the following fact.

\begin{bfhpg}[Fact]
  \label{Tfact}
  Let $T \xra{\tau} P \xra{\pi} M$ and $T' \xra{\tau'} P' \xra{\pi'}
  M'$ be complete projective resolutions. For every morphism
  $\mapdef{\a}{M}{M'}$ there exists a morphism $\widebar{\a}$ such
  that the right-hand square in the diagram
  \begin{equation*}
    \xymatrix{
      T \ar[r]^-{\tau} \ar[d]^-{\widehat{\a}} & P \ar[r]^-{\pi}
      \ar[d]^-{\widebar{\a}} & M \ar[d]^-{\a}\\
      T' \ar[r]^-{\tau'} & P' \ar[r]^-{\pi'} & M'}
  \end{equation*}
  is commutative up to homotopy. The morphism $\widebar{\a}$ is unique
  up to homotopy, and for every choice of $\widebar{\a}$ there exists
  a morphism $\widehat{\a}$, also unique up to homotopy, such that the
  left-hand square is commutative up to homotopy. Moreover, if $\tau'$
  and $\pi'$ are surjective, then $\widehat{\a}$ and $\widebar{\a}$
  can be chosen such that the diagram is commutative. Finally, if one
  has $M=M'$ and $\a$ is the identity map, then $\widebar{\a}$ and
  $\widehat{\a}$ are homotopy equivalences.
\end{bfhpg}

\begin{bfhpg}[Gorenstein projectivity]
  An $R$-module $G$ is called \emph{Gorenstein projective} if there
  exists a totally acyclic complex $T$ of projective $R$-modules with
  $\Co[0]{T} \is G$. In that case, the diagram $T \to \Thb{0}{T} \to
  G$ is a complete projective resolution, and for brevity we shall
  often say that $T$ is a complete projective resolution of $G$.

  The \emph{Gorenstein projective dimension} of an $R$-complex $M$,
  written $\Gpd{M}$, is the least integer $n$ such that there exists a
  complete projective resolution \eqref{cpltres} where $\tau_i$ is an
  isomorphism for all $i \ge n$. In particular, $\Gpd{M}$ is finite if
  and only if $M$ has a complete projective resolution. Notice that
  $\H{M}$ is bounded above if $\Gpd{M}$ is finite; indeed, there is an
  inequality
  \begin{equation}
    \label{eq:Gpd-sup}
    \Gpd{M} \ge \sup{\H{M}}.
  \end{equation}
  If $M$ is an $R$-complex of finite projective dimension, then there
  is a semi-projective resolution $P \qra M$ with $P$ bounded above,
  and then $0 \to P \to M$ is a complete projective resolution; in
  particular, $M$ has finite Gorenstein projective dimension.
\end{bfhpg}

\begin{bfhpg}[Tate homology]
  \label{Ttor}
  Let $M$ be an $\Rop$-complex with a complete projective resolution
  $T \to P \to M$. For an $R$-complex $N$, the \emph{Tate homology of
    $M$ with coefficients in $N$} is defined as
  \begin{equation*}
    \Ttor{i}{M}{N} =  \H[i]{\tp{T}{N}}.
  \end{equation*}
  It follows from \pgref{Tfact} that this definition is independent
  (up to isomorphism) of the choice of complete projective resolution;
  in particular, one has
  \begin{equation}
    \label{eq:Ttor=tor}
    \Ttor{i}{M}{N} \is \Tor{i}{M}{N}\quad\text{for } i > \Gpd[\Rop]{M}
    + \sup{N}.
  \end{equation}
  Note that $\Ttor{i}{M}{N}$ is a $\Bbbk$-module for every
  $i\in\ZZ$. Moreover, if $N$ is an \bi{R}{S}module, then each
  $\Ttor{i}{M}{N}$ is an $\Sop$-module.
\end{bfhpg}

Tate homology $\Ttor{*}{M}{N}$ vanishes if $M$ (or $N$) is a (bounded
above) complex of finite projective dimension; this is the content of
\prpref[]{vanTtor} and \lemref[]{f1} below.

The boundedness condition on $N$ in \lemref{f1} is a manifestation of
the fact that Tate homology $\Ttor{*}{M}{-}$ is not a functor from the
derived category $\D$. Indeed, every $R$-complex is isomorphic in $\D$
to a semi-projective complex, and for such a complex $P$ one has
$\Ttor{*}{M}{P}=0$ for every $\Rop$-complex $M$ of finite Gorenstein
projective dimension.

Notice, though, that if $M$ and $M'$ are isomorphic in $\D[\Rop]$ and
of finite Gorenstein projective dimension, then it follows from
\cite[1.4.P]{LLAHBF91} that every complete projective resolution $T
\to P \to M$ yields a complete resolution $T \to P \to M'$, so one has
an isomorphism $\Ttor{*}{M}{-} \is \Ttor{*}{M'}{-}$ of functors from
$\C$.

\begin{prp}
  \label{prp:vanTtor}
  Let $M$ be an $\Rop$-complex of finite Gorenstein projective
  dimension.  Among the conditions
  \begin{eqc}
  \item $\pd[\Rop]{M} < \infty$
  \item $\Ttor{i}{M}{-}=0$ for all $i\in\ZZ$
  \item $\Ttor{i}{M}{-}=0$ for some $i\in\ZZ$
  \item $\fd[\Rop]{M} < \infty$
  \end{eqc}
  the implications
  \eqclbl{i}$\!\implies\!$\eqclbl{ii}$\!\implies\!$\eqclbl{iii}$%
  \!\implies\!$\eqclbl{iv} hold.
\end{prp}

We recall from works of Jensen \prpcite[6]{CUJ70} and Raynaud and
Gruson \cite[II.\ thm.\ 3.2.6]{LGrMRn71} that if $R$ has finite
finitistic projective dimension---for example, $R$ is commutative
Noetherian of finite Krull dimension---then every flat $R$-module has
finite projective dimension, and it follows that the conditions
\eqclbl{i}--\eqclbl{iv} are equivalent.

\begin{prf*}
  If $\mapdef[\qra]{\pi}{P}{M}$ is a semi-projective resolution with
  $P$ bounded above, then $0 \to P \xra{\pi} M$ is a complete
  projective resolution, so one has $\Ttor{i}{M}{-}=0$ for all
  $i\in\ZZ$. Thus, \eqclbl{i} implies \eqclbl{ii}; the implication
  \proofofimp[]{ii}{iii} is trivial.

  Assume now that one has $\Ttor{i}{M}{-}=0$ for some $i\in\ZZ$. Let
  \mbox{$T \to P \to M$} be a complete projective resolution and set
  $G=\Co[i-1]{T}$. As the functor $\Tor{1}{G}{-} = \H[i]{\tp{T}{-}} =
  \Ttor{i}{M}{-}$ vanishes, the $\Rop$-module $G$ is flat. It follows
  that the module $\Co[j]{T} \is \Co[j]{P}$ is flat for every $j\ge
  \max\set{i-1,\Gpd[\Rop]{M}}$, whence $\fd[\Rop]{M}$ is finite.
\end{prf*}

\begin{rmk}
  \label{rmk:i1}
  Let $T \to P \to M$ be a complete projective resolution over
  $\Rop$. For every semi-projective resolution
  $\mapdef[\qra]{\pi'}{P'}{N}$ over $R$, application of the functor
  $\tp{T}{-}$ to the exact sequence $0 \to N \to \Cone{\pi'} \to
  \Shift{P'} \to 0$ yields a short exact sequence, as $T$ is a complex
  of projective $\Rop$-modules. The associated exact sequence in
  homology yields an isomorphism
  \begin{equation}
    \label{eq:Ttoraug}
    \H{\tp{T}{N}} \is \H{\tp{T}{\Cone{\pi'}}},
  \end{equation}
  as one has $\H{\tp{T}{P'}}=0$ because $P'$ is semi-flat. If $N$ is
  bounded above and of finite projective dimension, then one can
  assume that $P'$ and, therefore, $\Cone{\pi'}$ is bounded above, and
  then \lemcite[2.13]{CFH-06} yields
  $\H{\tp{T}{\Cone{\pi'}}}=0$. Thus, we record the following result.
\end{rmk}

\begin{lem}
  \label{lem:f1}
  Let $M$ be an $\Rop$-complex of finite Gorenstein projective
  dimension. For every bounded above $R$-complex $N$ of finite
  projective dimension, one has $\Ttor{i}{M}{N}=0$ for all
  $i\in\ZZ$.\qed
\end{lem}

\enlargethispage*{\baselineskip}

\begin{prp}
  \label{prp:Ttor-3N}
  Let $M$ be an $\Rop$-complex of finite Gorenstein projective
  dimension. For every exact sequence $0 \to N' \to N \to N'' \to 0$
  of $R$-complexes, there is an exact sequence of\, $\Bbbk$-modules
  \begin{equation*}
    \cdots \to \Ttor{i+1}{M}{N''} \to \Ttor{i}{M}{N'} \to
    \Ttor{i}{M}{N} \to  \Ttor{i}{M}{N''} \to \cdots.
  \end{equation*}
  Moreover, if the original exact sequence is one of complexes of
  \bi{R}{S}modules, then the derived exact sequence is one of
  $\Sop$-modules.
\end{prp}

\begin{prf*}
  Let $T \to P \to M$ be a complete projective resolution. The
  sequence
  \begin{equation*}
    0 \to \tp{T}{N'} \to \tp{T}{N} \to \tp{T}{N''} \to 0
  \end{equation*}
  is exact because $T$ is a complex of projective $\Rop$-modules. The
  associated exact sequence in homology is the desired one, and the
  statement about additional module structures is evident.
\end{prf*}

\begin{prp}
  \label{prp:Ttor-3M}
  Let $0 \to M' \to M \to M'' \to 0$ be an exact sequence of
  $\Rop$\nobreakdash-complexes of finite Gorenstein projective
  dimension. For every $R$-complex $N$ there is an exact sequence of\,
  $\Bbbk$-modules
  \begin{equation*}
    \cdots \to \Ttor{i+1}{M''}{N} \to \Ttor{i}{M'}{N} \to
    \Ttor{i}{M}{N} \to  \Ttor{i}{M''}{N} \to \cdots.
  \end{equation*}
  Moreover, if $N$ is a complex of \bi{R}{S}modules, then the derived
  exact sequence is one of $\Sop$-modules.
\end{prp}

\begin{prf*}
  By \prpcite[4.7]{OVl06} there is a commutative diagram
  \begin{equation*}
    \xymatrix{
      0 \ar[r] & T' \ar[r] \ar[d] & T \ar[r] \ar[d] & T'' \ar[r]
      \ar[d] & 0\\
      0 \ar[r] & P' \ar[r] \ar[d] & P \ar[r] \ar[d] & P'' \ar[r]
      \ar[d] & 0\\
      0 \ar[r] & M' \ar[r] & M \ar[r] & M'' \ar[r] & 0}
  \end{equation*}
  with exact rows, such that the columns are complete projective
  resolutions. The sequence $0 \to T' \to T \to T'' \to 0$ is
  degreewise split, so the sequence
  \begin{equation*}
    0 \to \tp{T'}{N} \to \tp{T}{N} \to \tp{T''}{N} \to 0
  \end{equation*}
  is exact, and the associated sequence in homology is the desired
  one. The statement about additional module structures is evident.
\end{prf*}

As with absolute homology, dimension shifting is a useful technique in
dealings with Tate homology.

\begin{lem}
  \label{lem:Ttor-dimshift}
  Let $M$ be an $\Rop$-complex of finite Gorenstein projective
  dimension and let $N$ be an $R$-complex. For every complete
  projective resolution $T \to P \to M$ and for every $m\in\ZZ$ there
  are isomorphisms
  \begin{equation*}
    \tag{a}
    \Ttor{i}{M}{N} \is \Ttor{i-m}{\Co[m]{T}}{N} \ \text{ for all $i\in\ZZ$.}
  \end{equation*}
  For every semi-projective resolution $L \qra N$ and for every
  integer $n \ge \sup{N}$ there are isomorphisms
  \begin{equation*}
    \tag{b}
    \Ttor{i}{M}{N} \is \Ttor{i-n}{M}{\Co[n]{L}} \ \text{ for all $i\in\ZZ$.}
  \end{equation*}
\end{lem}

\begin{prf*}
  (a): For every $m\in\ZZ$ the diagram $\Shift[-m]{T} \to
  \Shift[-m]{\Thb{m}{T}} \to \Co[m]{T}$ is a complete projective
  resolution. Hence one has
  \begin{align*}
    \Ttor{i-m}{\Co[m]{T}}{N} &= \H[i-m]{\tp{\Shiftp[-m]{T}}{N}}\\
    &= \H[i-m]{\Shift[-m]{\tpp{T}{N}}}\\
    &\is \H[i]{\tp{T}{N}}\\ & = \Ttor{i}{M}{N}.
  \end{align*}

  (b): We may assume that $N$ is bounded above; otherwise the
  statement is void.  For every $n \ge \sup{N}$ there is a
  quasi-isomorphism $\mapdef{\widetilde{\pi}}{\Tsa{n}{L}}{N}$. The
  acyclic complex $\Cone{\widetilde{\pi}}$ is bounded above, so
  $\tp{T}{\Cone{\widetilde{\pi}}}$ is acyclic by
  \lemcite[2.13]{CFH-06}. An application of \prpref{Ttor-3N} to the
  exact sequence $0 \to N \to \Cone{\widetilde{\pi}} \to
  \Shift{\Tsa{n}{L}} \to 0$ yields isomorphisms
  \begin{equation*}
    \Ttor{i}{M}{N} \is \Ttor{i}{M}{\Tsa{n}{L}} \ \text{ for all $i\in\ZZ$.}
  \end{equation*}
  The first complex in the exact sequence $0 \to \Tha{n-1}{L} \to
  \Tsa{n}{L} \to \Shift[n]{\Co[n]{L}} \to 0$ of $R$-complexes has
  finite projective dimension. Indeed, in the exact sequence $0 \to
  \Tha{n-1}{L} \to L \to \Thb{n}{L} \to 0$, the complexes $L$ and
  $\Thb{n}{L}$ are semi-projective, so $\Tha{n-1}{L}$ is
  semi-projective and, moreover, bounded above. Now apply \lemref{f1}
  and \prpref{Ttor-3N} to get
  \begin{equation*}
    \Ttor{i}{M}{\Tsa{n}{L}} \is \Ttor{i}{M}{\Shift[n]{\Co[n]{L}}} \
    \text{ for all $i\in\ZZ$.}
  \end{equation*}
  The desired isomorphisms follow from these last two displays.
\end{prf*}

\section{Pinched tensor product complexes}
\label{sec:box}

\noindent
We start by noticing that a very natural approach to the balancedness
problem for Tate homology fails.

\begin{exa}
  \label{exa:TT}
  Let $k$ be a field and consider the commutative ring
  $R=\poly{x,y}/(xy)$. The $R$-module $R/(x)$ is Gorenstein projective
  with complete resolution
  \begin{equation*}
    T = \cdots \xra{y} R \xra{x} R \xra{y} R \xra{x} \cdots,
  \end{equation*}
  where $\dif[i]{T}$ is multiplication by $x$ for $i$ odd and
  multiplication by $y$ for $i$ even. As multiplication by $y$ on
  $R/(x)$ is injective, it is immediate from the definition of Tate
  homology, see \pgref{Ttor}, that one has $\Ttor{i}{R/(x)}{R/(x)}=0$
  for $i$ even.

  The complex $\tp{T}{T}$, however, has non-vanishing homology in even
  degrees. Indeed, for each $n\in\ZZ$ the module $\tpp{T}{T}_n$ is
  free with basis $(e_{i,n-i})_{i\in\ZZ}$. The differential is given
  by
  \begin{equation*}
    \dif[n]{\tp{T}{T}}(e_{i,n-i}) = 
    \begin{cases}
      xe_{i-1,n-i} - ye_{i,n-i-1}&\text{$n$ odd and $i$ odd}\\
      ye_{i-1,n-i} + xe_{i,n-i-1}&\text{$n$ odd and $i$ even}\\
      xe_{i-1,n-i} - xe_{i,n-i-1}&\text{$n$ even and $i$ odd}\\
      ye_{i-1,n-i} + ye_{i,n-i-1}&\text{$n$ even and $i$ even.}
    \end{cases}
  \end{equation*}
  For $n$ even, the element $xe_{0,n}$ is a cycle and clearly not a
  boundary. Indeed, since $R$ is graded, the complex $\tp{T}{T}$ has
  an internal grading, and the differential is of degree $1$ with
  respect to this grading.  Suppose that $xe_{0,n}$ is a boundary.
  Since it is an element of internal degree 1, a preimage
  $\sum_{i\in\mathbb Z}\alpha_{i,n+1-i} e_{i,n+1-i}$ of $xe_{0,n}$
  under $\dif{\tp{T}{T}}$ may be assumed homogeneous of internal
  degree zero.  That is, we may assume that $\alpha_{i,n+1-i}$ is in
  $k$ for all $i$.  Let $i_0$ and $i_1$ be, respectively, the least
  and the largest integer $i$ with $\alpha_{i,n+1-i}\ne 0$.  With
  respect to the basis $(e_{i,n-i})_{i\in\ZZ}$, the element
  $b=\dif{\tp{T}{T}}(\sum_{i\in\mathbb Z}\alpha_{i,n+1-i}
  e_{i,n+1-i})$ is nonzero in coordinate $(i_0-1,n+1-i_0)$, which
  implies $i_0=1$.  Similarly, $b$ is nonzero in coordinate
  $(i_1,n-i_1)$, which implies $i_1=0$.  Thus one has $i_0>i_1$, a
  contradiction.
\end{exa}

The isomorphism \eqref{Ttoraug} shows, nevertheless, that one can
compute Tate homology from a tensor product of acyclic complexes. This
motivates the next construction; see also the comments before the
proof of \thmref{btp-Ttor}.

\begin{con}
  \label{con:btp}
  Let $T$ be an $\Rop$-complex and let $A$ be an $R$-complex. Consider
  the graded $\Bbbk$-module $\btp{T}{A}$ defined by:
  \begin{equation*}
    \btpp{T}{A}_n =
    \begin{cases}
      \tpp{\Thb{0}{T}}{\Thb{0}{A}}_n & \text{ for $n \ge 0$}\\
      \tpp{\Tha{-1}{T}}{\Shift{\Thap{-1}{A}}}_n & \text{ for $n \le
        -1$.}
    \end{cases}
  \end{equation*}
  It is elementary to verify that one has
  \begin{equation*}
    (\tp{\dif[0]{T}}{(\sigma\dif[0]{A})}) \circ
    \dif[1]{\tp{\Thb{0}{T}}{\Thb{0}{A}}} =0 =
    \dif[-1]{\tp{\Tha{-1}{T}}{\Shift{\Thap{-1}{A}}}} \circ
    (\tp{\dif[0]{T}}{(\sigma\dif[0]{A})}),
  \end{equation*}
  where $\sigma$ denotes the canonical map $A \to \Shift{A}$.  Thus,
  $\dif{\btp{T}{A}}$ defined by
  \begin{equation*}
    \dif[n]{\btp{T}{A}} =
    \begin{cases}
      \dif[n]{\tp{\Thb{0}{T}}{\Thb{0}{A}}} & \text{ for $n \ge 1$}\\
      \tp{\dif[0]{T}}{(\sigma\dif[0]{A})} &  \text{ for $n =0$}\\
      \dif[n]{\tp{\Tha{-1}{T}}{\Shift{\Thap{-1}{A}}}} & \text{ for $n
        \le -1$}
    \end{cases}
  \end{equation*}
  is a differential on $\btp{T}{A}$. We refer to this $\Bbbk$-complex
  as the \emph{pinched tensor product} of $T$ and $A$.

  For morphisms $\mapdef{\a}{T}{T'}$ of $\Rop$-complexes and
  $\mapdef{\b}{A}{A'}$ of $R$-complexes, it is elementary to verify
  that the assignment $x\otimes y \mapsto \a(x) \otimes \b(y)$ defines
  a morphism of $\Bbbk$-complexes
  \begin{equation*}
    \dmapdef{\btp{\a}{\b}}{\btp{T}{A}}{\btp{T'}{A'}}.
  \end{equation*}
\end{con}

\begin{rmk}
  \label{rmk:btp}
  For every $\Rop$-complex $T$ and every $R$-complex $A$ there are
  equalities of $\Bbbk$-complexes,
  \begin{align}
    \label{eq:btp1}
    \Thb{0}{\btpp{T}{A}} &= \tp{\Thb{0}{T}}{\Thb{0}{A}}
    \qand\\
    \label{eq:btp2}
    \Tha{-1}{\btpp{T}{A}} &= \tp{\Tha{-1}{T}}{\Shift{\Thap{-1}{A}}}.
  \end{align}
  If $T$ is a complex of \bi{R'}{R}modules, and $A$ is a complex of
  \bi{R}{S}modules, then $\btp{T}{A}$ is a complex of
  \bi{R'}{S}modules.
\end{rmk}

The proof of the next proposition is standard, and we omit it.

\begin{prp}
  The pinched tensor product defined in {\rm \conref[]{btp}} yields a
  functor
  \begin{equation*}
    \dmapdef{\btp{-}{-}}{\C[\bicat{R'}{R}] \times
      \C[\bicat{R}{S}]}{\C[\bicat{R'}{S}]};
  \end{equation*}
  in particular, it yields a functor $\C[\Rop]\times\C \to
  \C[\Bbbk]$. Moreover, it is $\Bbbk$-bilinear and right exact in each
  variable.  \qed
\end{prp}

\begin{thm}
  \label{thm:btp-Ttor}
  Let $M$ be an $\Rop$-complex with a complete projective resolution
  $T \to P \to M$. Let $A$ be an acyclic $R$-complex and set $N =
  \Co[0]{A}$. For every $i \in \ZZ$ there is an isomorphism of\,
  $\Bbbk$-modules
  \begin{equation*}
    \H[i]{\btp{T}{A}}\cong \Ttor{i}{M}{N}.
  \end{equation*}
  If $A$ is a complex of \bi{R}{S}modules, then the isomorphism is one
  of $\Sop$-modules.
\end{thm}

Before we proceed with the proof, we point out that if $N$ is an
$R$-module, and $A$ is the acyclic complex $0 \to N \xra{=} N \to 0$
with $N$ in degrees $0$ and $-1$, then one has $\btp{T}{A} =
\tp{T}{N}$.

\begin{prf*}
  By definition one has $\Ttor{i}{M}{N} = \H[i]{\tp{T}{N}}$, so the
  goal is to establish an isomorphism between $\H{\btp{T}{A}}$ and
  $\H{\tp{T}{N}}$.  The quasi-isomorphisms
  \begin{equation*}
    \pi\colon \Thb{0}{A} \qra N \qqand \epsilon\colon N \qra
    \Shift{\Thap{-1}{A}},
  \end{equation*}
  with $\epsilon_0\pi_0 = \sigma\dif[0]{A}$, induce quasi-isomorphisms
  \begin{equation*}
    \Thb{0}{\btpp{T}{A}} \qra \tp{\Thb{0}{T}}{N} \qand
    \tp{\Tha{-1}{T}}{N} \qra \Tha{-1}{\btpp{T}{A}};
  \end{equation*}
  see \eqref{btp1}, \eqref{btp2}, and \prpcite[2.14]{CFH-06}. It
  follows that there are isomorphisms $\H[i]{\btp{T}{A}} \is
  \H[i]{\tp{T}{N}}$ for all $i\in\ZZ\setminus\set{0,-1}$. To establish
  the isomorphism in the remaining two degrees, consider the following
  diagram with exact columns.
  \begin{equation*}
    \xymatrixrowsep{1.35pc}
    \xymatrixcolsep{1.35pc}
    \xymatrix{
      & 0 \ar[d] & & \\
      & \tp{T_0}{\Bo[0]{A}} \ar[d] & 0 \ar[d] & 0 \ar[d] \\
      \btpp{T}{A}_1 \ar[r] \ar[d]_{\tp[]{T_{1}}{\pi_0}} &
      \btpp{T}{A}_0 \ar[dr] \ar[d]_{\tp[]{T_0}{\pi_0}} &
      \tpp{T}{N}_{-1} \ar[r] \ar[d]^{\tp[]{T_{-1}}{\epsilon_0}} &
      \tpp{T}{N}_{-2} \ar[d]^{\tp[]{T_{-2}}{\epsilon_0}} \\
      \tpp{T}{N}_1 \ar[r] \ar[d] &
      \tpp{T}{N}_0 \ar[d] \ar[ur] &
      \btpp{T}{A}_{-1}  \ar[r] \ar[d]^-{\tp[]{T_{-1}}{\dif[-1]{A}}} &
      \btpp{T}{A}_{-2}\\
      0 & 0 & \tp{T_{-1}}{\Bo[-2]{A}}
      \ar[d]\\ & & 0 }
  \end{equation*}
  The identity $\epsilon_0\pi_0=\sigma\dif[0]{A}$ shows that the
  twisted square is commutative. That the other two squares are
  commutative follows by functoriality of the tensor product.

  To see that the homomorphism $\tp{T_0}{\pi_0}$ induces the desired
  isomorphism in homology, $\H[0]{\btp{T}{A}} \is \H[0]{\tp{T}{N}}$,
  notice first that it maps boundaries to boundaries, and that for $x$
  in $\Cy[0]{\btp{T}{A}}$ one has
  \begin{equation*}
    \tpp{T_{-1}}{\epsilon_0}\circ\dif[0]{\tp{T}{N}}\circ\tpp{T_0}{\pi_0}(x)=0
  \end{equation*}
  by commutativity of the twisted square. As $\tp{T_{-1}}{\epsilon_0}$
  is injective, it follows that $\tpp{T_0}{\pi_0}(x)$ is in
  $\Cy[0]{\tp{T}{N}}$, so there is a well-defined homomorphism
  \begin{equation*}
    \H{\tp{T_0}{\pi_0}}\colon \H[0]{\btp{T}{A}} \to \H[0]{\tp{T}{N}}.
  \end{equation*}
  It is immediate from the surjectivity of $\tp{T_0}{\pi_0}$ and
  commutativity of the twisted square that the homomorphism
  $\H{\tp{T_0}{\pi_0}}$ is surjective. To see that it is injective,
  let $x$ be an element in $\Cy[0]{\btp{T}{A}}$ and assume that there
  is a $y$ in $\tpp{T}{N}_1$ such that $\tpp{T_0}{\pi_0}(x) =
  \dif[1]{\tp{T}{N}}(y)$. Choose an element $z$ in $\tp{T_1}{A_0}
  \subset \btpp{T}{A}_1$ such that $\tpp{T_1}{\pi_0}(z)=y$. Then the
  element $x- \dif[1]{\btp{T}{A}}(z)$ in $\btpp{T}{A}_0$ maps to $0$
  under $\tp{T_0}{\pi_0}$, so it belongs to $\tp{T_0}{\Bo[0]{A}}$. Let
  $w$ in $\tp{T_0}{A_1} \subset \btpp{T}{A}_1$ be a preimage of $x-
  \dif[1]{\btp{T}{A}}(z)$. Then one has
  \begin{equation*}
    \dif[1]{\btp{T}{A}}(w)=(\tp{T_0}{\dif[1]{A}})(w)=x-\dif[1]{\btp{T}{A}}(z),
  \end{equation*}
  and so $x$ is a boundary: $\dif[1]{\btp{T}{A}}(w+z)=x$. Thus,
  $\H{\tp{T_0}{\pi_0}}$ is an isomorphism.

  Similarly, for $i=-1$, it is evident that $\tp{T_{-1}}{\epsilon_0}$
  maps cycles to cycles. Let $x$ be a boundary in $\tpp{T}{N}_{-1}$,
  and choose a preimage $y$ of $x$ in $\tpp{T}{N}_0$.  By surjectivity
  of $\tp{T_0}{\pi_0}$, this $y$ has a preimage $z$ in
  $\btpp{T}{A}_0$, and by commutativity of the twisted square one has
  $\dif[0]{\btp{T}{A}}(z) = \tpp{T_{-1}}{\epsilon_0}(x)$. Thus,
  $\tp{T_{-1}}{\epsilon_0}$ maps boundaries to boundaries, whence it
  induces a homomorphism
  \begin{equation*}
    \H{\tp{T_{-1}}{\epsilon_0}}\colon \H[-1]{\tp{T}{N}} \to \H[-1]{\btp{T}{A}}.
  \end{equation*}
  It follows immediately from the injectivity of
  $\tp{T_{-1}}{\epsilon_0}$ and commutativity of the twisted square
  that $\H{\tp{T_{-1}}{\epsilon_0}}$ is injective. To see that it is
  surjective, let $x$ be an element in $\Cy[-1]{\btp{T}{A}}$.  Then,
  in particular, one has
  \begin{equation*}
    0=\tpp{T_{-1}}{\dif[0]{\Shift{\Thap{-1}{A}}}}(x)=-\tpp{T_{-1}}{\dif[-1]{A}}(x).
  \end{equation*}
  Therefore, $x$ is in $\tp{T_{-1}}{\Cy[-1]{A}} =
  \Im(T_{-1}\otimes\epsilon_0)$, and it follows by injectivity of
  $\tp{T_{-2}}{\epsilon_0}$ that the preimage of $x$ is a cycle in
  $\tp{T_{-1}}{N}$. Thus, $\H{\tp{T_{-1}}{\epsilon_0}}$ is surjective
  and hence an isomorphism.

  The claim about $\Sop$-module structures is immediate from
  \conref{btp}.
\end{prf*}

\begin{prp}
  \label{prp:btp-comm}
  Let $T$ be an $\Rop$-complex and let $A$ be an $R$-complex. The map
  \begin{equation*}
    \varpi\colon \btp{T}{A} \lra \btp[\Rop]{A}{T}
  \end{equation*}
  given by
  \begin{alignat*}{2}
    \varpi_n (t \otimes a) &= (-1)^{\dgr{t}\dgr{a}}a\otimes
    t & &\text{for $n \ge 0$}\\
    \varpi_n(t \otimes \sigma(a)) &= (-1)^{(\dgr{t}+1)(\dgr{a}+1)}
    a\otimes\sigma(t) &\quad &\text{for $n \le -1$}
  \end{alignat*}
  is an isomorphism of\, $\Bbbk$-complexes.

  Moreover, if\, $T$ is a complex of \bi{R'}{R}modules and $A$ is a
  complex of \bi{R}{S}mo-dules, then $\varpi$ is an isomorphism of
  complexes of \bi{R'}{S}modules.
\end{prp}

\begin{prf*}
  The map $\varpi$ is clearly an isomorphism of graded
  $\Bbbk$-modules, and it is straightforward to verify that it
  commutes with the differentials. The assertions about additional
  module structures are immediate from \conref{btp}.
\end{prf*}

If $M$ is an $\Rop$-module of finite Gorenstein projective dimension
and $N$ is an $R$-module of finite Gorenstein projective dimension,
then one could also define Tate homology of the pair $(M,N)$ in terms
of the complete projective resolution of $N$. Do the two definitions
agree; that is, is Tate homology balanced? This is tantamount to
asking if one has $\Ttor{*}{M}{N} \is \Ttor[{\Rop}]{*}{N}{M}$. Iacob
\cite{AIc07} gave a positive answer for modules over commutative
Noetherian Gorenstein rings. The next theorem settles the question
over any associative ring.

\begin{thm}
  \label{thm:Ttor-bal}
  Let $M$ be an $\Rop$-complex and let $N$ be an $R$-complex, both of
  which are both bounded above and of finite Gorenstein projective
  dimension. For every $i\in\ZZ$ there is an isomorphism of\,
  $\Bbbk$-modules:
  \begin{equation*}
    \Ttor{i}{M}{N} \is \Ttor[{\Rop}]{i}{N}{M}.
  \end{equation*}
\end{thm}

\begin{prf*}
  Choose complete projective resolutions $T \to P \to M$ and $T' \to
  P' \to N$. Set $m = \max\set{\sup{M},\Gpd[\Rop]{M}}$ and $n =
  \max\set{\sup{N},\Gpd{N}}$.  The modules $\Co[m]{P} \is \Co[m]{T}$
  and $\Co[n]{P'}\is \Co[n]{T'}$ are Gorenstein projective with
  complete projective resolutions
  \begin{equation*}
    \Shift[-m]{T} \to \Shift[-m]{\Thb{m}{P}} \to \Co[m]{P} \qqand
    \Shift[-n]{T'} \to \Shift[-n]{\Thb{n}{P'}} \to \Co[n]{P'}.
  \end{equation*}
  \lemref{Ttor-dimshift}, \thmref{btp-Ttor}, and \prpref{btp-comm} now
  conspire to yield the desired isomorphism,
  \begin{align*}
    \Ttor{i}{M}{N} 
    &\is \Ttor{i-m-n}{\Co[m]{P}}{\Co[n]{P'}} \\
    &\is \H[i-m-n]{\btp{\Shiftp[-m]{T}}{\Shift[-n]{T'}}} \\ &\is
    \H[i-n-m]{\btp[R^\circ]{\Shiftp[-n]{T'}}{\Shift[-m]{T}}} \\
    &\is \Ttor[\Rop]{i-n-m}{\Co[n]{P'}}{\Co[m]{P}} \\
    &\is \Ttor[\Rop]{i}{N}{M}. \qedhere
  \end{align*}
\end{prf*}

\begin{rmk}
  In \cite{AIc07} Iacob considers a variation of Tate homology based
  on complete flat resolutions. The proof of \thmref{btp-Ttor}
  applies, mutatis mutandis, to show that also these homology groups
  can be computed from a pinched tensor product. From a result
  parallel to \lemref{Ttor-dimshift} it, therefore, follows that also
  this version of Tate homology is balanced.
\end{rmk}

\section{Pinched Hom complexes and Tate cohomology}
\label{sec:tatecohomology}

\noindent
Tate cohomology was studied in detail by Veliche~\cite{OVl06}; we
recall the definition.

\begin{ipg}
  \label{TC}
  Let $M$ be an $R$-complex with a complete projective resolution $T
  \to P \to M$. For an $R$-complex $N$, the \emph{Tate cohomology of
    $M$ with coefficients in $N$} is defined~as
  \begin{equation*}
    \Text{i}{M}{N} =  \H[-i]{\Hom{T}{N}}.
  \end{equation*}
  This definition is independent (up to isomorphism) of the choice of
  complete resolution; cf.~\pgref{Tfact}. In particular, one has
  \begin{equation}
    \label{eq:Text=ext}
    \Text{i}{M}{N} \is \Ext{i}{M}{N}\quad\text{for } i > \Gpd{M} - \inf{N}.
  \end{equation}

  Note that $\Text{i}{M}{N}$ is a $\Bbbk$-module for every
  $i\in\ZZ$. Moreover, if $N$ is an \bi{R}{S}module, then each
  $\Text{i}{M}{N}$ is an $\Sop$-module.
\end{ipg}

The parallels of \prpref[]{vanTtor}--\prpref[]{Ttor-3M} are
established in \seccite[4]{OVl06}. The proof of \lemref{Text-dimshift}
is similar to the proof of \lemref{Ttor-dimshift}. It uses
\lemcite[2.4]{CFH-06} and the following fact, which follows from an
argument similar to the one given in \rmkref{i1}.

\begin{lem}
  \label{lem:f2}
  \protect\pushQED{\qed} Let $M$ be an $R$-complex of finite
  Gorenstein projective dimension.  For every bounded below
  $R$-complex $N$ of finite injective dimension, one has
  $\Text{i}{M}{N}=0$ for all $i\in\ZZ$. \qed
\end{lem}

\begin{lem}
  \protect\pushQED{\qed}
  \label{lem:Text-dimshift}
  Let $M$ be an $R$-complex of finite Gorenstein projective dimension
  and let $N$ be an $R$-complex. For every complete projective
  resolution $T \to P \to M$ and for every integer $m \in \ZZ$ there
  are isomorphisms
  \begin{equation*}
    \tag{a}
    \Text{i}{M}{N} \is \Text{i-m}{\Co[m]{T}}{N} \ \text{ for all $i\in\ZZ$.}
  \end{equation*}
  For every semi-injective resolution $N \qra I$ and for every integer
  $n \ge -\inf{N}$ there are isomorphisms
  \begin{equation}
    \tag{b}
    \Text{i}{M}{N} \is \Text{i-n}{M}{\Cy[-n]{I}} \ \text{ for all
      $i\in\ZZ$.} \qedhere
  \end{equation}
\end{lem}

\begin{con}
  \label{con:pHom}
  Let $T$ and $A$ be $R$-complexes. Consider the graded $\Bbbk$-module
  $\pHom{T}{A}$ defined by:
  \begin{equation*}
    \pHom{T}{A}_n =
    \begin{cases}
      \Hom{\Tha{-1}{T}}{\Shift[-1]{\Thbp{1}{A}}}_n & \text{ for $n \ge 1$} \\
      \Hom{\Thb{0}{T}}{\Tha{0}{A}}_n & \text{ for $n \le 0$.}
    \end{cases}
  \end{equation*}
  It is elementary to verify that one has
  \begin{equation*}
    \Hom{\dif[0]{T}}{\dif[1]{A}\varsigma} \circ
    \dif[2]{\Hom{\Tha{-1}{T}}{\Shift[-1]{\Thbp{1}{A}}}} = 0 =
    \dif[0]{\Hom{\Thb{0}{T}}{\Tha{0}{A}}} \circ
    \Hom{\dif[0]{T}}{\dif[1]{A}\varsigma},
  \end{equation*}
  where $\varsigma$ denotes the canonical map $\Shift[-1]{A} \to
  A$. Thus, $\dif[n]{\pHom{T}{A}}$ defined by
  \begin{equation*}
    \dif[n]{\pHom{T}{A}} =
    \begin{cases}
      \dif[n]{\Hom{\Tha{-1}{T}}{\Shift[-1]{\Thbp{1}{A}}}} & \text{ for
        $n \ge 2$}\\
      \Hom{\dif[0]{T}}{\dif[1]{A}\varsigma} &  \text{ for $n=1$}\\
      \dif[n]{\Hom{\Thb{0}{T}}{\Tha{0}{A}}} & \text{ for $n \le 0$}
    \end{cases}
  \end{equation*}
  is a differential on $\pHom{T}{A}$. We refer to this $\Bbbk$-complex
  as the \emph{pinched Hom} of $T$ and~$A$.

  For morphisms $\mapdef{\a}{T}{T'}$ and $\mapdef{\b}{A}{A'}$ of
  $R$-complexes it is elementary to verify that the assignment $\f
  \mapsto \b\f\a$ defines a morphism of $\Bbbk$-complexes
  \begin{equation*}
    \dmapdef{\pHom{\a}{\b}}{\pHom{T'}{A}}{\pHom{T}{A'}}.
  \end{equation*}
\end{con}

\begin{rmk}
  \label{rmk:pHom}
  For all $R$-complexes $T$ and $A$ there are equalities of complexes
  \begin{align}
    \label{eq:pHom1}
    \Thb{1}{\pHom{T}{A}} &= \Hom{\Tha{-1}{T}}{\Shift[-1]{\Thbp{1}{A}}}
    \qand\\
    \label{eq:pHom2}
    \Tha{0}{\pHom{T}{A}} &= \Hom{\Thb{0}{T}}{\Tha{0}{A}}.
  \end{align}
  If $T$ is a complex of \bi{R}{R'}modules, and $A$ is a complex of
  \bi{R}{S}modules, then $\pHom{T}{A}$ is a complex of
  \bi{R'}{S}modules.
\end{rmk}

The proof of the next proposition is standard and omitted.

\begin{prp}
  The pinched Hom defined in {\rm \conref[]{pHom}} yields a functor
  \begin{equation*}
    \dmapdef{\pHom{-}{-}}{\C[\bicat{R}{R'}]^\mathrm{op} \times
      \C[\bicat{R}{S}]}{\C[\bicat{R'}{S}]};
  \end{equation*}
  in particular, it yields a functor $\C^\mathrm{op}\times\C \to
  \C[\Bbbk]$. Moreover, it is $\Bbbk$-bilinear and left exact in each
  variable.  \qed
\end{prp}

\begin{thm}
  \label{thm:pHom-Text1}
  Let $M$ be an $R$-complex with a complete projective resolution
  $T\to P\to M$. Let $A$ be an acyclic $R$-complex and set $N =
  \Cy[0]{A}$.  For every $i \in \ZZ$ there is an isomorphism of\,
  $\Bbbk$-modules
  \begin{equation*}
    \H[-i]{\pHom{T}{A}}\cong \Text{i}{M}{N}.
  \end{equation*}
  If $A$ is a complex of \bi{R}{S}modules, then the isomorphism is one
  of $\Sop$-modules.
\end{thm}

Notice that if $N$ is an $R$-module and $A$ is the acyclic complex $0
\to N \xra{=} N \to 0$ with $N$ in degrees $1$ and $0$, then one has
$\pHom{T}{A} = \Hom{T}{N}$.

\begin{prf*}
  The quasi-isomorphisms
  \begin{equation*}
    \pi\colon\Shift[-1]{\Thbp{1}{A}} \qra N \qqand \epsilon\colon N
    \qra \Tha{0}{A},
  \end{equation*}
  with $\epsilon_0\pi_0=\dif[1]{A}\varsigma$, yield quasi-isomorphisms
  \begin{equation*}
    \Thb{1}{\pHom{T}{A}} \qra \Hom{\Tha{-1}{T}}{N} \qand
    \Hom{\Thb{0}{T}}{N} \qra \Tha{0}{\pHom{T}{A}};
  \end{equation*}
  see \eqref{pHom1}, \eqref{pHom2}, and \prpcite[2.6]{CFH-06}. It
  follows that there are isomorphisms $\H[i]{\pHom{T}{A}} \is
  \H[i]{\Hom{T}{N}}$ for all $i\in\ZZ\setminus\set{1,0}$. To establish
  the desired isomorphism for $i\in\set{0,1}$, consider the following
  diagram with exact columns.
  \begin{equation*}
    \xymatrixrowsep{1.35pc}
    \xymatrixcolsep{1.35pc}
    \xymatrix{
      & 0 \ar[d] & \\
      &
      \mspace{-30mu}\Hom{T_{-1}}{\Cy[1]{A}}\mspace{-30mu}
      \ar[d] &  0\ar[d] & 0 \ar[d] \\
      \pHom{T}{A}_2 \ar[r] \ar[d]_-{\Hom[]{T_{-2}}{\pi_0}} &
      \pHom{T}{A}_1 \ar[d]_-{\Hom[]{T_{-1}}{\pi_0}} \ar[dr]  & \Hom{T}{N}_{0} \ar[r]
      \ar[d]^-{\Hom[]{T_0}{\epsilon_0}} & \Hom{T}{N}_{-1}
      \ar[d]^-{\Hom[]{T_1}{\epsilon_0}}\\
      \Hom{T}{N}_2 \ar[d] \ar[r]&
      \Hom{T}{N}_1 \ar[d] \ar[ur] & \pHom{T}{A}_{0} \ar[r]
      \ar[d]^-{\Hom{T_0}{\dif[0]{A}}}& \pHom{T}{A}_{-1}\\^
      0 & 0 & \mspace{-16mu}\Hom{T_0}{\Bo[-1]{A}}\mspace{-16mu} \ar[d]\\
      & & 0}
  \end{equation*}
  The identity $\epsilon_0\pi_0 = \dif[1]{A}\varsigma$ ensures that
  the twisted square is commutative; also the other two squares are
  commutative by functoriality of the Hom functor.

  To see that $\Hom{T_{-1}}{\pi_0}$ induces an isomorphism from
  $\H[1]{\pHom{T}{A}}$ to $\H[1]{\Hom{T}{N}}$, notice first that it
  maps boundaries to boundaries by commutativity of the left-hand
  square. For a cycle $\zeta$ in $\Cy[1]{\pHom{T}{A}}$ one has
  \begin{equation*}
    \big(\Hom{T_{0}}{\epsilon_0}\circ\dif[1]{\Hom{T}{N}} \circ
    \Hom{T_{-1}}{\pi_0}\big)(\zeta)=0,
  \end{equation*}
  by commutativity of the twisted square. As $\Hom{T_{0}}{\epsilon_0}$
  is injective, it follows that $\Hom{T_{-1}}{\pi_0}(\zeta)$ is in
  $\Cy[1]{\Hom{T}{N}}$, so there is a well-defined homomorphism
  \begin{equation*}
    \H{\Hom{T_{-1}}{\pi_0}}\colon \H[1]{\pHom{T}{A}} \to
    \H[1]{\Hom{T}{N}}.
  \end{equation*}
  It is immediate by surjectivity of $\Hom{T_{-1}}{\pi_0}$ and
  commutativity of the twisted square that $\H{\Hom{T_{-1}}{\pi_0}}$
  is surjective. To see that it is also injective, let $\zeta$ be a
  cycle in $\pHom{T}{A}_1$ and assume that one has
  $\Hom{T_{-1}}{\pi_0}(\zeta) = \dif[2]{\Hom{T}{N}}(\alpha)$ for some
  element $\alpha$ in $\Hom{T}{N}_2 = \Hom{T_{-2}}{N}$. For any
  preimage $\xi$ of $\alpha$ in $\Hom{T_{-2}}{A_1} \subset
  \pHom{T}{A}_2$, the element $\zeta- \dif[2]{\pHom{T}{A}}(\xi)$ in
  $\pHom{T}{A}_1$ maps to $0$ under $\Hom{T_{-1}}{\pi_0}$, so it
  belongs to $\Hom{T_{-1}}{\Cy[1]{A}}$. As $T_{-1}$ is projective and
  $A$ is acyclic, there exists a homomorphism $\psi$ in
  $\Hom{T_{-1}}{A_2} \subset \pHom{T}{A}_2$ such that one has
  \begin{equation*}
    \zeta- \dif[2]{\pHom{T}{A}}(\xi) =
    \dif[2]{A}\psi =
    -\dif[2]{\pHom{T}{A}}(\psi).
  \end{equation*}
  It follows that $\zeta$ is a boundary, $\zeta =
  \dif[2]{\pHom{T}{A}}(\xi - \psi)$, whence $\H{\Hom{T_{-1}}{\pi_0}}$
  is an isomorphism.

  From the commutativity of the right-hand square, it follows that
  $\Hom{T_0}{\epsilon_0}$ maps cycles to cycles. Let $\beta$ be a
  boundary in $\Hom{T}{N}_0$ and choose a preimage $\alpha$ of $\beta$
  in $\Hom{T}{N}_1$.  By surjectivity of $\Hom{T_{-1}}{\pi_0}$ this
  $\alpha$ has a preimage $\alpha'$ in $\pHom{T}{A}_1$, and by
  commutativity of the twisted square one has
  $\Hom{T_{0}}{\epsilon_0}(\beta) =
  \dif[1]{\pHom{T}{A}}(\alpha')$. Thus, $\Hom{T_{0}}{\epsilon_0}$ maps
  boundaries to boundaries, whence it it induces a homomorphism
  \begin{equation*}
    \H{\Hom{T_{0}}{\epsilon_0}}\colon \H[0]{\Hom{T}{N}} \to \H[0]{\pHom{T}{A}}.
  \end{equation*}
  It follows immediately from injectivity of $\Hom{T_{0}}{\epsilon_0}$
  and commutativity of the twisted square that
  $\H{\Hom{T_0}{\epsilon_0}}$ is injective. To see that it is
  surjective, let $\zeta$ be a cycle in $\pHom{T}{A}_0$; one then has
  $0 = \dif[0]{\Tha{0}{A}}\zeta = \Hom{T_0}{\dif[0]{A}}(\zeta)$.  By
  exactness of the second column from the right, it now follows that
  $\zeta$ is in the image of $\Hom{T_0}{\epsilon_0}$, and by
  injectivity of $\Hom{T_1}{\epsilon_0}$ it follows that the preimage
  of $\zeta$ is a cycle in $\Hom{T}{N}_{0}$. Thus,
  $\H{\Hom{T_{0}}{\epsilon_0}}$ is an isomorphism.

  The claim about $\Sop$-module structures is immediate from
  \conref{pHom}.
\end{prf*}

The next result is a pinched version of Hom-tensor adjunction.

\begin{prp}
  \label{prp:padj}
  Let $T$ be an $R$-complex, let $A$ be a complex of \bi{S}{R}modules,
  and let $B$ be an $S$-module. The map
  \begin{equation*}
    \dmapdef{\varrho}{\Hom[S]{\btp[\Rop]{T}{A}}{B}}{\pHom{T}{\Hom[S]{A}{B}}}
  \end{equation*}
  given by
  \begin{align*}
    \varrho_n(\psi)(t)(a) =
    \begin{cases}
      \psi(t \otimes \sigma(a)) & \text{for
        $n\ge 1$}\\
      \psi(t \otimes a) &\text{for $n\le 0$} \\
    \end{cases}
  \end{align*}
  is an isomorphism of\, $\Bbbk$-complexes.

  Moreover, if $T$ is a complex of \bi{R}{R'}modules, and $B$ is an
  \bi{S}{S'}module, then $\varrho$ is an isomorphism of complexes of
  \bi{R'}{S'}modules.
\end{prp}

\begin{prf*} 
  For $n \ge 1$ one has
  \begin{align*}
    \Hom[S]{\btp[\Rop]{T}{A}}{B}_n 
    &= \Hom[S]{(\btp[\Rop]{T}{A})_{-n}}{B}\\
    &=
    \Hom[S]{\textstyle\dfinsum_{i=-n}^{-1}\tp[\Rop]{T_{i}}{\Shiftp{A}_{-n-i}}}{B}
  \end{align*}
  and
  \begin{align*}
    \textstyle\dfinsum_{i=-n}^{-1}
    \Hom{T_{i}}{\Hom[S]{A_{-n-i-1}}{B}} 
    &=\textstyle\dfinsum_{i=-n}^{-1}
    \Hom{T_{i}}{\Hom[S]{A}{B}_{i+n+1}} \\
    &=\textstyle\dfinsum_{i=-n}^{-1} \Hom{T_{i}}{\Shiftp[-1]{\Hom[S]{A}{B}}_{i+n}} \\
    &=\pHom{T}{\Hom[S]{A}{B}}_n.
  \end{align*}
  The map $\varrho_n$ given by $\varrho_n(\psi)(t)(a) = \psi(t \otimes
  \sigma(a))$, for $t\in T_{i}$ and $a\in A_{-n-i-1}$ is, up to
  $\sigma$, just the Hom-tensor adjunction isomorphism of
  modules. Thus, $\varrho_n$ is an isomorphism of
  $\Bbbk$-modules. Moreover, still for $n\ge 1$, one has
  \begin{multline*}
    \varrho_n\big(\dif[n+1]{\Hom[S]{\btp[\Rop]{T}{A}}{B}}(\psi)\big)(t)(a)\\
    \begin{aligned}
      & = \varrho_n\big(-(-1)^{n+1}\psi\dif[-n]{\btp[\Rop]{T}{A}}\big)(t)(a)\\
      & = (-1)^{n}\psi\dif[-n]{\btp[\Rop]{T}{A}}(t\otimes \sigma(a))\\
      & = (-1)^{n}\psi\big( \dif{T}(t)\otimes \sigma(a)
      +(-1)^{\dgr{t}}t\otimes
      \dif{\Shift{A}}(\sigma(a)\big)\\
      & = (-1)^{n}\psi\big( \dif{T}(t)\otimes \sigma(a)
      -(-1)^{\dgr{t}}t\otimes
      \sigma(\dif{A}(a))\big)\\
      & = (-1)^{n+\dgr{t}+1}\psi(t\otimes \sigma(\dif{A}(a))) +
      (-1)^{n}\psi(\dif{T}(t)\otimes \sigma(a))
    \end{aligned}
  \end{multline*}
  and
  \begin{multline*}
    \big(\dif[n+1]{\pHom{T}{\Hom[S]{A}{B}}}\varrho_{n+1}(\psi)\big)(t)(a) \\
    \begin{aligned}
      & = \big(\dif{\Shift[-1]{\Hom[S]{A}{B}}}\varrho_{n+1}(\psi) -
      (-1)^{n+1} \varrho_{n+1}(\psi)\dif{T}\big)(t)(a)\\
      & =
      (-1)^{\dgr{\varrho_{n+1}(\psi)(t)}}\varrho_{n+1}(\psi)(t)(\dif{A}(a))
      -(-1)^{n+1}\varrho_{n+1}(\psi)(\dif{T}(t))(a)\\
      & = (-1)^{n+1 +\dgr{t}}\psi(t\otimes\sigma(\dif{A}(a)))+
      (-1)^{n}\psi(\dif{T}(t) \otimes \sigma(a)).
    \end{aligned}
  \end{multline*}
  That is, the identity
  $\dif[n+1]{\pHom{T}{\Hom[S]{A}{B}}}\varrho_{n+1} =
  \varrho_n\dif[n+1]{\Hom[S]{\btp[\Rop]{T}{A}}{B}}$ holds for $n \ge
  1$.

  By \eqref{btp1} and \eqref{pHom1} there are equalities of
  $\Bbbk$-complexes
  \begin{align*}
    \Tha{0}{\Hom[S]{\btp[\Rop]{T}{A}}{B}} &=
    \Hom[S]{\tp[\Rop]{\Thb{0}{T}}{\Thb{0}{A}}}{B}\qand\\
    \Tha{0}{\pHom{T}{\Hom[S]{A}{B}}}
    &
    = \Hom{\Thb{0}{T}}{\Hom[S]{\Thb{0}{A}}{B}}.
  \end{align*}
  Thus, for $n \le 0$ the map $\varrho_n$ is the degree $n$ component
  of the Hom-tensor adjunction isomorphism
  $\Hom[S]{\tp[\Rop]{\Thb{0}{T}}{\Thb{0}{A}}}{B} \xra{\is}
  \Hom{\Thb{0}{T}}{\Hom[S]{\Thb{0}{A}}{B}}$.
  
  To prove that $\varrho$ is an isomorphism of $\Bbbk$-complexes, it
  remains to verify the identity
  $\dif[1]{\pHom{T}{\Hom[S]{A}{B}}}\varrho_{1} =
  \varrho_0\dif[1]{\Hom[S]{\btp[\Rop]{T}{A}}{B}}$. For $t\in T_0$ and
  $a\in A_0$ one has
  \begin{align*}
    \varrho_0\big(\dif[1]{\Hom[S]{\btp[\Rop]{T}{A}}{B}}(\psi)\big)(t)(a)
    &= \varrho_0\big(\psi\dif[0]{\btp[\Rop]{T}{A}}\big)(t)(a)\\
    &= \psi(\dif{T}(t) \otimes \sigma\dif{A}(a))
  \end{align*}
  and
  \begin{align*}
    \big(\dif[1]{\pHom{T}{\Hom[S]{A}{B}}}\varrho_{1}(\psi)\big)(t)(a)
    &= \big(\dif[1]{\Hom[S]{A}{B}}\varsigma\varrho_1(\psi)(\dif{T}(t))\big)(a)\\
    &= \varrho_1(\psi)(\dif{T}(t))(\dif{A}(a))\\
    &= \psi(\dif{T}(t) \otimes \sigma\dif{A}(a))
  \end{align*}
  as required. Finally, the statements about extra module structures
  are evident in view of the remarks made in \rmkref[]{btp} and
  \rmkref[]{pHom}.
\end{prf*}

\begin{prp}
  \label{prp:pHom-btp}
  Assume that $R$ is commutative.  Let $M$ be an $R$-complex with a
  complete projective resolution $T \to P \to M$ and let $N$ be a
  Gorenstein projective $R$-module with complete projective resolution
  $T'$.  For every projective $R$-module $Q$ and every $i \in \ZZ$
  there is an isomorphism of $R$-modules
  \begin{equation*}
    \H[-i]{\Hom{\btp{T}{T'}}{Q}}\cong \Text{i}{M}{\Hom{N}{Q}}.
  \end{equation*}
\end{prp}

\begin{prf*}
  The $R$-complex $\Hom{T'}{Q}$ is acyclic, and $\Hom{N}{Q}$ is the
  kernel of the differential in degree $0$. The assertion now follows
  from \prpref{padj} and \thmref{pHom-Text1}.
\end{prf*}

\begin{cor}
  \label{cor:pp=p}
  Assume that $R$ is commutative.  Let $M$ and $N$ be Gorenstein
  projective $R$-modules with complete projective resolutions $T$ and
  $T'$, respectively. If one has $\Ttor{i}{M}{N}=0$ for all $i\in\ZZ$,
  then the complex $\btp{T}{T'}$ of projective $R$-modules is acyclic,
  and the following conditions are equivalent.
  \begin{eqc}
  \item The $R$-complex $\btp{T}{T'}$ is totally acyclic.
  \item For every projective $R$-module~$Q$ one has
    $\Text{i}{M}{\Hom{N}{Q}}=0$ for all $i\in\ZZ$.
  \end{eqc}
  When these conditions hold, the $R$-module $\tp{M}{N}$ is Gorenstein
  projective with complete projective resolution $\btp{T}{T'}$.
\end{cor}

\begin{prf*}
  By construction the complex $\btp{T}{T'}$ consists of projective
  $R$-modules, and one has $\Co[0]{\btp{T}{T'}} \is \tp{M}{N}$. The
  assumption that the Tate homology $\Ttor{*}{M}{N}$ vanishes implies
  that $\btp{T}{T'}$ is acyclic; see \thmref{btp-Ttor}. The
  equivalence of \eqclbl{i} and \eqclbl{ii} now follows from
  \prpref{pHom-btp}, and the last assertion is then evident.
\end{prf*}

\section{Tate cohomology is balanced}

\noindent
For $R$-modules $M$ and $N$, a potentially different approach to Tate
cohomology $\Text{*}{M}{N}$ uses a resolution of the second argument
$N$. The resulting theory, which is parallel to the one developed
in~\cite{LLAAMr02,ROB86,OVl06}, was outlined by Asadollahi and
Salarian in \cite{JAsSSl06}. In this section we use the pinched
complexes to show that when both approaches apply, they yield the same
cohomology theory.

\begin{bfhpg}[Complete injective resolutions]
  A complex $U$ of injective $R$\nobreakdash-mo\-dules is called
  \emph{totally acyclic} if it is acyclic, and the complex
  $\Hom{J}{U}$ is acyclic for every injective $R$-module $J$.

  A \emph{complete injective resolution} of an $R$-complex $N$ is a
  diagram
  \begin{equation}
    \label{eq:cpltcores}
    N \xra{\iota} I \xra{\upsilon} U,
  \end{equation}
  where $\iota$ is a semi-injective resolution, $U$ is a totally
  acyclic complex of injective $R$-modules, and $\upsilon_i$ is an
  isomorphism for $i \ll 0$.
\end{bfhpg}

\begin{bfhpg}[Gorenstein injectivity]
  An $R$-module $E$ is called \emph{Gorenstein injective} if there
  exists a totally acyclic complex $U$ of injective $R$-modules with
  $\Cy[0]{U} \is E$. In that case, the diagram $E \to \Tha{0}{U} \to
  U$ is a complete injective resolution, and for brevity we shall
  often say that $U$ is a complete injective resolution of $E$.

  The \emph{Gorenstein injective dimension} of an $R$-complex $N$,
  written $\Gid{N}$, is the least integer $n$ such that there exists a
  complete injective resolution \eqref{cpltcores} where $\upsilon_i$
  is an isomorphism for all $i \le -n$. In particular, $\Gid{N}$ is
  finite if and only if $N$ has a complete injective
  resolution. Notice that $\H{N}$ is bounded below if $\Gid{N}$ is
  finite; indeed, there is an inequality
  \begin{equation}
    \label{eq:Gid-inf}
    \Gid{N} \ge -\inf{\H{N}}.
  \end{equation}
  If $N$ is an $R$-complex of finite injective dimension, then there
  is a semi-injective resolution $N \qra I$ with $I$ bounded below,
  and then $N \to I \to 0$ is a complete injective resolution; in
  particular, $N$ has finite Gorenstein injective dimension.
\end{bfhpg}

\begin{prp}
  \label{prp:pHom-Text2}
  Let $N$ be an $R$-complex with a complete injective resolution $N
  \to I \to U$. Let $A$ be an acyclic $R$-complex and set $M =
  \Co[0]{A}$.  For every $i \in \ZZ$ there is an isomorphism of\,
  $\Bbbk$-modules
  \begin{equation*}
    \H[i]{\pHom{A}{U}} \cong \H[i]{\Hom{M}{U}}.
  \end{equation*}
  If $A$ is a complex of \bi{R}{S}modules, then the isomorphism is one
  of $S$-modules.
\end{prp}

Notice that if $M$ is an $R$-module and $A$ is the acyclic complex $0
\to M \xra{=} M \to 0$ with $M$ in degrees $0$ and $-1$, then one has
$\pHom{A}{U} \is \Hom{M}{U}$.

\begin{prf*}
  The quasi-isomorphisms
  \begin{equation*}
    \pi\colon{\Thb{0}{A}} \qra M \qqand  \epsilon\colon M \qra
    \Shift{\Thap{-1}{A}},
  \end{equation*}
  with $\epsilon_0\pi_0=\sigma\dif[0]{A}$, induce quasi-isomorphisms
  \begin{align*}
    \Hom{M}{\Tha{0}{U}} &\qra \Tha{0}{\pHom{A}{U}} \qand\\
    \Hom{\Shift{\Tha{-1}{A}}}{\Thb{1}{U}} &\qra \Hom{M}{\Thb{1}{U}};
  \end{align*}
  see \eqref{pHom1} and \prpcite[2.7]{CFH-06}. There is an equality of
  graded $\Bbbk$-modules
  \begin{equation*}
    \Hom{\Shift{\Tha{-1}{A}}}{\Thb{1}{U}} = \Thb{1}{\pHom{A}{U}},
  \end{equation*}
  and one has $-\dif{\Hom{\Shift{\Tha{-1}{A}}}{\Thb{1}{U}}} =
  \dif{\Thb{1}{\pHom{A}{U}}}$. It follows that there are isomorphisms
  $\H[i]{\pHom{A}{U}} \is \H[i]{\Hom{M}{U}}$ for $i\ne 1,0$. To
  establish the desired isomorphism in the remaining two degrees,
  consider the following diagram with exact columns.
  \begin{equation*}
    \xymatrixrowsep{1.35pc}
    \xymatrixcolsep{1.23pc}
    \xymatrix{
      & 0 \ar[d] &  \\
      & \mspace{-48mu}\Hom{\Bo[-2]{A}}{U_1}\mspace{-48mu}
      \ar[d]_-{\Hom[]{\dif[-1]{A}}{U_1}} & 0 \ar[d] & 0 \ar[d]\\
      \pHom{A}{U}_2 \ar[r] \ar[d]_-{\!-\Hom[]{\epsilon_0}{U_2}}&
      \pHom{A}{U}_1 \ar[d]_-{\Hom[]{\epsilon_0}{U_1}} \ar[dr]  &
      \Hom{M}{U}_{0} \ar[r]
      \ar[d]^-{\Hom[]{\pi_0}{U_0}} & \Hom{M}{U}_{-1}
      \ar[d]^-{\Hom[]{\pi_0}{U_{-1}}}\! \\
      \Hom{M}{U}_2 \ar[d] \ar[r]&
      \Hom{M}{U}_1 \ar[d] \ar[ur] & 
      \pHom{A}{U}_{0} \ar[r] \ar[d] & \pHom{A}{U}_{-1} \\
      0 & 0 & \mspace{-36mu}\Hom{\Bo[0]{A}}{U_0}\mspace{-36mu} \ar[d] \\ 
      & & 0 }
  \end{equation*}
  The identity $\epsilon_0\pi_0=\sigma\dif[0]{A}$ ensures that the
  twisted square is commutative; also the other two squares are
  commutative by standard properties of the Hom functor.

  To see that $\Hom{\epsilon_0}{U_1}$ and $\Hom{\pi_0}{U_0}$ induce
  isomorphisms in homology, one proceeds as in the proof of
  \thmref{pHom-Text1}.
\end{prf*}

If $M$ is a Gorenstein projective $R$-module with complete projective
resolution $T$, and $N$ is a Gorenstein injective $R$-module with
complete injective resolution $U$, then \thmref{pHom-Text1} and
\prpref{pHom-Text2} yield
\begin{equation*}
  \Text{i}{M}{N} \is \H[-i]{\pHom{T}{U}} \is \H[-i]{\Hom{M}{U}}.
\end{equation*}
That is, the Tate cohomology of $M$ with coefficients in $N$ can be
computed via a complete injective resolution of $N$. What follows is a
balancedness statement that shows that for appropriately bounded
complexes---for modules in particular---one can unambiguously extend
the notion of Tate cohomology $\Text{*}{M}{N}$ to the situation where
$N$ has a complete injective resolution; see \dfnref{Text}.

\begin{thm}
  \label{thm:Text-bal}
  Let $M$ be a bounded above $R$-complex with a complete projective
  resolution and let $N$ be a bounded below $R$-complex with a
  complete injective resolution $N \to I \to U$. For every $i\in\ZZ$
  there is an isomorphism
  \begin{equation*}
    \Text{i}{M}{N} \is \H[-i]{\Hom{M}{U}}.
  \end{equation*}
\end{thm}

\begin{prf*}
  Set $n = \sup\set{-\inf{N},\Gid{N}}$; then the module $\Cy[-n]{I}
  \is \Cy[-n]{U}$ is Gorenstein injective with complete injective
  resolution $\Cy[-n]{I} \to \Shift[n]{\Tha{-n}{I}} \to
  \Shift[n]{U}$. Further, set $m = \Gpd{M}$ and let $T \to P \to M$ be
  a complete projective resolution; then the module $\Co[m]{P} \is
  \Co[m]{T}$ is Gorenstein projective with complete projective
  resolution $\Shift[-m]{T}\to \Shift[-m]{\Thb{m}{P}} \to
  \Co[m]{P}$. In the next chain of isomorphisms, the first one follows
  from \lemref{Text-dimshift}, the second and third follow from
  \thmref{pHom-Text1} and \prpref{pHom-Text2}, and the last one
  follows by dimension shifting.
  \begin{align*}
    \Text{i}{M}{N} &\is \Text{i-m-n}{\Co[m]{P}}{\Cy[-n]{I}} \\
    &\is \H[m+n-i]{\pHom{\Shift[-m]{T}}{\Shift[n]{U}}} \\
    & \is \H[m+n-i]{\Hom{\Co[m]{P}}{\Shift[n]{U}}}\\
    & \is \H[m-i]{\Hom{\Co[m]{P}}{U}}.
  \end{align*}
  Finally, an argument parallel to the one for
  \lemref{Ttor-dimshift}(b) yields isomorphisms
  \begin{equation*}
    \H[-i]{\Hom{M}{U}} \is \H[-i]{\Hom{\Tsa{m}{P}}{U}} \is 
    \H[m-i]{\Hom{\Co[m]{P}}{U}};
  \end{equation*}
  this time it is \lemcite[2.5]{CFH-06} that needs to be invoked.
\end{prf*}

\begin{dfn}
  \label{dfn:Text}
  Let $N$ be a bounded below $R$-complex with a complete injective
  resolution $N \to I \to U$. For every bounded above $R$-complex $M$,
  the \emph{Tate cohomology of $M$ with coefficients in $N$} is given
  by
  \begin{equation*}
    \Text{i}{M}{N} =  \H[-i]{\Hom{M}{U}}.
  \end{equation*}
\end{dfn}

\begin{rmk}
  A fact parallel to \pgref{Tfact} guarantees that the definition
  above is independent (up to isomorphism) of the choice of complete
  resolution; in particular, one has the following parallel of
  \eqref{Text=ext},
  \begin{equation}
    \label{eq:Text=ext2}
    \Text{i}{M}{N} \is \Ext{i}{M}{N}\quad\text{for } i > \Gid{N} + \sup{M}.
  \end{equation} 
  Other standard results similar to
  \prpref[]{vanTtor}--\prpref[]{Ttor-3M} are established in
  \cite{JAsSSl06}. In that paper, the notation
  $\smash{\operatorname{\widehat{ext}}}_{R}^{*{^{\phantom{|}\mspace{-6mu}}}}(M,N)$
  is used for the cohomology defined in \dfnref[]{Text}, and it is
  shown to agree with the notion from \cite{LLAAMr02,ROB86,OVl06}, see
  \pgref{TC}, over commutative Noetherian local Gorenstein~rings.

  More generally, for a module $N$ with a complete injective
  resolution, Nucinkis' \cite{BEN98} notion of I-complete cohomology
  agrees with Tate cohomology as defined in
  \dfnref[]{Text}. Similarly, for a module $M$ with a complete
  projective resolution, the P-complete cohomology of Benson and
  Carlson \cite{DFBJFC92}, Vogel/Goichot \cite{FGc92}, and Mislin
  \cite{GMs94} agrees with Tate cohomology in the sense of
  \pgref{TC}. Nucinkis proves \thmcite[5.2, 6.6, 7.9]{BEN98} that P-
  and I-complete cohomology agree over rings where every module has a
  complete projective resolution and a complete injective resolution.
\end{rmk}

The next result establishes a pinched version of the Hom swap
isomorphism. It is proved in the same fashion as \prpref{padj}.

\begin{prp}
  \label{prp:pswap}
  Let $T$ be an $R$-complex, let $B$ be an $\Sop$-module, and let $U$
  be a complex of \bi{R}{S}modules.  The map
  \begin{equation*}
    \dmapdef{\vartheta}{\Hom[\Sop]{B}{\pHom{T}{U}}}{\pHom{T}{\Hom[\Sop]{B}{U}}}
  \end{equation*}
  given by
  \begin{equation*}
    \vartheta_n(\psi)(t)(b) = \psi(b)(t)
  \end{equation*}
  is an isomorphism of\, $\Bbbk$-complexes.

  Moreover, if $T$ is a complex of \bi{R}{R'}modules, and $B$ is an
  \bi{S'}{S}module, then $\vartheta$ is an isomorphism of complexes of
  \bi{R'}{S'}modules.\qed
\end{prp}

\begin{prp}
  \label{prp:pHom-btp1}
  Assume that $R$ is commutative. Let $M$ be an $R$-complex with a
  complete projective resolution $T \to P \to M$ and let $N$ be a
  Gorenstein injective $R$-module with complete injective resolution
  $U$. For every injective $R$-module $J$ and every $i\in\ZZ$ there is
  an isomorphism of $R$-modules
  \begin{equation*}
    \H[-i]{\Hom{J}{\pHom{T}{U}}}\cong \Text{i}{M}{\Hom{J}{N}}.
  \end{equation*}
\end{prp}

\begin{prf*}
  The complex $\Hom{J}{U}$ is acyclic and $\Hom{J}{N}$ is the kernel
  of the differential in degree $0$. The assertion now follows from
  \prpref{pswap} and \thmref{pHom-Text1}.
\end{prf*}

\begin{cor}
  \label{cor:pi=i}
  Assume that $R$ is commutative. Let $M$ be a Gorenstein projective
  $R$-module with complete projective resolution $T$ and let $N$ be a
  Gorenstein injective $R$-module with complete injective resolution
  $U$. If one has $\Text{i}{M}{N}=0$ for all $i\in\ZZ$, then the
  complex $\pHom{T}{U}$ of injective $R$-modules is acyclic, and the
  following conditions are equivalent.
  \begin{eqc}
  \item The $R$-complex $\pHom{T}{U}$ is totally acyclic.
  \item For every injective $R$-module $J$ one has
    $\Text{i}{M}{\Hom{J}{N}}=0$ for all $i\in\ZZ$.
  \end{eqc}
  When these conditions hold, the $R$-module $\Hom{M}{N}$ is
  Gorenstein injective with complete injective resolution
  $\pHom{T}{U}$.
\end{cor}

\begin{prf*}
  By construction the complex $\pHom{T}{U}$ consists of injective
  $R$\nobreakdash-modules, and one has $\Cy[0]{\pHom{T}{U}} \is
  \Hom{M}{N}$. The assumption that the Tate cohomology
  $\Text{*}{M}{N}$ vanishes implies that $\pHom{T}{U}$ is acyclic; see
  \thmref{pHom-Text1}. The equivalence of \eqclbl{i} and \eqclbl{ii}
  now follows from \prpref{pHom-btp1}, and the last assertion is then
  evident.
\end{prf*}

\section{Local algebra}

\noindent
Throughout this section $R$ denotes a commutative Noetherian local
ring with maximal ideal $\m$. Recall that every projective $R$-module
is free. An acyclic complex $T$ of finitely generated free
$R$\nobreakdash-modules is totally acyclic if and only if $\Hom{T}{R}$
is acyclic. For an $R$-module $M$ we use the standard notation $M^*$
for the dual module $\Hom{M}{R}$. A finitely generated $R$-module $G$
is Gorenstein projective if and only if one has
\begin{equation*}
  G \is G^{**} \qand \Ext{i}{G}{R} = 0 = \Ext{i}{G^*}{R} \ \text{ for all
    $i \ge 1$,}
\end{equation*}
see \cite{LLAAMr02}, and following \emph{op.\,cit.}\ we use the term
\emph{totally reflexive} for such modules.

A complex $F$ of finitely generated free $R$-modules is called
\emph{minimal} if one has $\partial(F) \subseteq \m F$; see
\seccite[8]{LLAAMr02}. A complete projective resolution $T \to P \to
M$ is called minimal if $T$ and $P$ are minimal complexes of finitely
generated free $R$\nobreakdash-modules. By \thmcite[8.4]{LLAAMr02}
every finitely generated $R$-module $M$ of finite Gorenstein
projective dimension has a minimal complete projective resolution
$T\to P \to M$, and it is unique up to isomorphism. The invariants
$\widehat\beta_n(M) = \rnk[R]{T_n}$ are called the \emph{stable Betti
  numbers} of $M$; for $n \ge \Gpd{M}$ they agree with usual Betti
numbers.

\begin{thm}
  \label{thm:TR}
  Let $M$ and $N$ be totally reflexive $R$-modules with complete
  projective resolutions $T$ and $T'$, respectively. If one has
  $\Ttor{i}{M}{N}=0$ for all $i\in\ZZ$, then the complex $\btp{T}{T'}$
  of finitely generated free $R$-modules is acyclic with
  $\Co[0]{\btp{T}{T'}} \is \tp{M}{N}$, and the following conditions
  are equivalent.
  \begin{eqc}
  \item The $R$-complex $\btp{T}{T'}$ is totally acyclic.
  \item One has $\Text{i}{M}{N^*}=0$ for all $i\in\ZZ$.
  \item The $R$-module $\tp{M}{N}$ is totally reflexive.
  \end{eqc}
  When these conditions hold, $\btp{T}{T'}$ is a complete projective
  resolution of $\tp{M\!}{\!N}$. It is minimal if and only if $T$ and
  $T'$ are minimal; in particular, one has
  \begin{equation*}
    \widehat\beta_i(\tp{M}{N}) = 
    \begin{cases}
      \sum_{0\le j \le i} \widehat\beta_j(M)\widehat\beta_{i-j}(N),& 
      \text{ for } i\ge 0\\[.75ex]
      \sum_{i\le j < 0} \widehat\beta_{j}(M)\widehat\beta_{i-j-1}(N),& 
      \text{ for } i < 0
    \end{cases}
  \end{equation*}
\end{thm}

\begin{prf*}
  By \conref{btp} the complex $\btp{T}{T'}$ consists of finitely
  generated free $R$-modules, and the assumption that Tate homology
  $\Ttor{*}{M}{N}$ vanishes implies that $\btp{T}{T'}$ is acyclic; see
  \thmref{btp-Ttor}.  To prove equivalence of the three conditions it
  suffices, in view of \corref{pp=p}, to prove the implication
  $\proofofimp[]{iii}{i}$. Assume that $\Co[0]{\btp{T}{T'}}=\tp{M}{N}$
  is totally reflexive. It follows immediately that the syzygies of
  $\tp{M}{N}$, i.e.\ $\Co[i]{\btp{T}{T'}}$ for $i \ge 1$ are totally
  reflexive as well. For $i \le -1$ it follows that
  $\Co[i]{\btp{T}{T'}}$ has finite Gorenstein projective dimension.
  The Krull dimension $d$ of $R$ is an upper bound for the Gorenstein
  projective dimension of any $R$-module, so $\Co[i]{\btp{T}{T'}}$ is
  totally reflexive as it is the $d$th syzygy of
  $\Co[i-d]{\btp{T}{T'}}$; see \thmcite[3.1]{CFH-06}. Thus, each
  module $\Co[i]{\btp{T}{T'}}$ is totally reflexive, and then
  $\btp{T}{T'}$ is totally acyclic by \lemcite[2.4]{LLAAMr02}.

  The assertions about minimality follow immediately from
  \conref{btp}, and so does the equality of stable Betti numbers.
\end{prf*}

\begin{cor}
  \label{cor:TR-Gor}
  Let $R$ be Gorenstein and let $M$ and $N$ be totally reflexive
  $R$\nobreakdash-modu\-les with (minimal) complete projective
  resolutions $T$ and $T'$, respectively. If one has
  $\Ttor{i}{M}{N}=0$ for all $i\in\ZZ$, then $\tp{M}{N}$ is totally
  reflexive with (minimal) complete resolution $\btp{T}{T'}$.
\end{cor}

\begin{prf*}
  As $R$ is Gorenstein, every acyclic complex of projective modules is
  totally acyclic; see \lemcite[2.4]{LLAAMr02}.
\end{prf*}

For modules $M$ and $N$ of finite Gorenstein projective dimension,
vanishing of Tate homology $\Ttor{*}{M}{N}$ yields information about
the complex $\Ltp{M}{N}$ that encodes the absolute homology
$\Tor{*}{M}{N}$; we pursue this line of investigation in
\cite{LWCDAJb}.  This paper we close with an interpretation of the
Tate homology modules $\Ttor 0{M}{N}$ and $\Ttor{-1}{M}{N}$ in terms
of a natural homomorphism.

\begin{prp}
  \label{prp:tev-Ttor}
  Let $M$ and $N$ be finitely generated $R$-modules. If $M$ is totally
  reflexive, then there is an exact sequence of $R$-modules
  \begin{equation*}
    0\to\Ttor 0{M}N \to M\otimes_R N \xra{\theta_{MN}} \Hom {M^*}N \to \Ttor
    {-1}{M}N\to 0,
  \end{equation*}
  where $\theta_{MN}$ is the natural homomorphism given by $x\otimes
  y\longmapsto [\f \mapsto \f(x)y]$.
\end{prp}

\begin{prf*}
  Let $T \to P \to M$ be a minimal complete projective resolution.
  The natural map $\mapdef{\theta_{FN}}{\tp{F}{N}}{\Hom{F^*}{N}}$ is
  an isomorphism for every finitely generated free $R$-module $F$.
  Let $\mapdef{\widetilde{\theta}}{\tp{T_0}{N}}{\Hom{T^*_{-1}}{N}}$ be
  the homomorphism given by $t\otimes y \longmapsto [\f \mapsto
  \f(\dif{T}(t))y]$, then the following diagram is commutative.
  \begin{equation*}
    \xymatrixrowsep{1.5pc}
    \xymatrixcolsep{1.5pc}
    \xymatrix{
      \cdots \ar[r]& \tp{T_1}{N} \ar[r] \ar[d]_=  &
      \tp{T_0}{N} \ar[r] \ar[d]_= &
      \tp{T_{-1}}N \ar[r] \ar[d]^-{\theta}_-{\is} & \tp{T_{-2}}N
      \ar[r] \ar[d]^-{\theta}_-{\is} & \cdots\\
      \cdots \ar[r]& \tp{T_1}{N} \ar[r] &
      \tp{T_0}{N} \ar[r]^-{\widetilde\theta} &
      \Hom{T^*_{-1}}N \ar[r] & \Hom{T_{-2}^*}{N} \ar[r]&
      \cdots}
  \end{equation*}
  Thus, Tate homology $\Ttor{*}{M}{N}$ can be computed from the bottom
  complex. Let $\mapdef{\pi}{T_0}{M}$ and
  $\mapdef{\epsilon}{M}{T_{-1}}$ be the natural homomorphisms with
  $\epsilon\pi = \dif[0]{T}$, and consider the commutative diagram
  \begin{equation*}
    \xymatrixrowsep{1.5pc}
    \xymatrixcolsep{1.5pc}
    \xymatrix{
      \cdots \ar[r]& \tp{T_1}{N} \ar[d] \ar[r] &
      \tp{T_0}{N} \ar[r]^-{\widetilde\theta} \ar[d]^{\tp[]{\pi}{N}}&
      \Hom{T_{-1}^*}N \ar[r] & \Hom{T_{-2}^*}{N} \ar[r]&
      \cdots\\
      & 0 \ar[r] &
      \tp MN \ar[d] \ar[r]^-\theta & \Hom {M^*}N
      \ar[u]_{\Hom[]{\Hom[]{\epsilon}{R}}{N}} \ar[r] & 0 \ar[u]\\
      &  & 0 & 0 \ar[u]& }
  \end{equation*}
  A straightforward diagram chase shows that the homomorphisms
  $\tp{\pi}{N}$ and ${\Hom{\Hom{\epsilon}{R}}{N}}$ induce isomorphisms
  in homology.
\end{prf*}

The next statement is proved similarly; see
\lemcite[5.8.(3)]{LLAAMr02}.

\begin{prp}
  \label{prp:hev-Ttor}
  Let $M$ and $N$ be finitely generated $R$-modules. If $M$ is totally
  reflexive, then there is an exact sequence of $R$-modules
  \begin{equation*}
    0\to\Text{-1}{M}{N} \to M^*\otimes_R N \xra{\nu_{MN}} 
    \Hom{M}{N} \to \Text{0}{M}{N}  \to 0,
  \end{equation*}
  where $\nu_{MN}$ is the natural homomorphism given by $\f\otimes y
  \longmapsto [x \mapsto \f(x)y]$.\qed
\end{prp}

\section*{Acknowledgment}

We thank the anonymous referee for several corrections and valuable
comments that helped improve the text.

\bibliographystyle{amsplain}

  \providecommand{\arxiv}[2][AC]{\mbox{\href{http://arxiv.org/abs/#2}{\sf
  arXiv:#2 [math.#1]}}}
  \providecommand{\oldarxiv}[2][AC]{\mbox{\href{http://arxiv.org/abs/math/#2}{\sf
  arXiv:math/#2
  [math.#1]}}}\providecommand{\MR}[1]{\mbox{\href{http://www.ams.org/mathscinet-getitem?mr=#1}{#1}}}
  \renewcommand{\MR}[1]{\mbox{\href{http://www.ams.org/mathscinet-getitem?mr=#1}{#1}}}
\providecommand{\bysame}{\leavevmode\hbox to3em{\hrulefill}\thinspace}
\providecommand{\MR}{\relax\ifhmode\unskip\space\fi MR }
\providecommand{\MRhref}[2]{%
  \href{http://www.ams.org/mathscinet-getitem?mr=#1}{#2}
}
\providecommand{\href}[2]{#2}

\end{document}